\documentclass[reqno,a4paper,12pt]{amsart}
\usepackage{latexsym, amsmath, amssymb, amsthm, a4}
\usepackage{graphicx}
\usepackage{subcaption}

\usepackage{amsfonts, bbold, bbm, marvosym}

\usepackage{mathrsfs}

\usepackage[mathscr]{eucal}
\font\sc=rsfs10 at 12pt
\renewcommand{\aa}{\overset{{}_{{}_\circ}}{\mathrm{a}}}

\renewcommand{\a}{\alpha}
\newcommand{\be}{\beta}
\newcommand{\g}{\gamma}
\newcommand{\G}{\Gamma}
\newcommand{\de}{\delta}
\newcommand{\D}{\Delta}
\newcommand{\e}{\epsilon}

\newcommand{\y}{\eta}

\newcommand{\io}{\iota}
\newcommand{\ka}{\kappa}

\newcommand{\la}{\lambda}

\newcommand{\m}{\mu}
\newcommand{\n}{\nu}
\newcommand{\x}{\xi}

\newcommand{\ro}{\rho}
\newcommand{\rob}{\pmb{\pmb{\rho}}}

\newcommand{\Si}{\Sigma}

\newcommand{\hi}{\chi}

\newcommand{\psb}{\pmb{\psi}}
\newcommand{\Psib}{\pmb{\Psi}}
\newcommand{\om}{\omega}
\newcommand{\Om}{\Omega}
\newcommand{\Omb}{\pmb{\Omega}}

\newcommand{\C}{{\mathbb C}}
\newcommand{\R}{{\mathbb R}}

\newcommand{\N}{{\mathbb{N}}}

\newcommand{\ab}{{\mathbf a}}

\newcommand{\bb}{{\mathbf b}}

\newcommand{\eb}{{\mathbf e}}

\newcommand{\kb}{{\mathbf k}}

\newcommand{\nb}{{\mathbf n}}

\newcommand{\qb}{{\mathbf q}}
\newcommand{\rb}{{\mathbf r}}

\newcommand{\vb}{{\mathbf v}}
\newcommand{\wb}{{\mathbf w}}

\newcommand{\Cb}{{\mathbf C}}

\newcommand{\Gb}{{\mathbf G}}
\newcommand{\Hb}{{\mathbf H}}
\newcommand{\Ib}{{\mathbf I}}

\newcommand{\Kb}{{\mathbf K}}

\newcommand{\Mb}{{\mathbf M}}
\newcommand{\Nb}{{\mathbf N}}

\newcommand{\Rb}{{\mathbf R}}
\newcommand{\Sbb}{{\mathbf S}}
\newcommand{\Tb}{{\mathbf T}}

\newcommand{\Vb}{{\mathbf V}}
\newcommand{\Wb}{{\mathbf W}}
\newcommand{\Xb}{{\mathbf X}}
\newcommand{\Yb}{{\mathbf Y}}
\newcommand{\Zb}{{\mathbf Z}}

\newcommand{\aF}{\mathfrak a}
\newcommand{\ap}{\pmb{\aF}}
\newcommand{\AF}{\mathfrak A}
\newcommand{\bFb}{\mathfrak b}
\newcommand{\BFB}{\mathfrak B}

\newcommand{\fF}{\mathfrak f}

\newcommand{\kF}{\mathfrak k}

\newcommand{\rF}{\mathfrak r}
\newcommand{\RF}{\mathfrak R}
\newcommand{\ssF}{\mathfrak s}

\newcommand{\WF}{\mathfrak W}

\newcommand{\XF}{\mathfrak X}

\newcommand{\UF}{\mathfrak U}

\newcommand{\Ac}{{\mathcal A}}

\newcommand{\Dc}{{\mathcal D}}
\newcommand{\Ec}{{\mathcal E}}
\newcommand{\Fc}{{\mathcal F}}
\newcommand{\Gc}{{\mathcal G}}

\newcommand{\Lc}{{\mathcal L}}
\newcommand{\Lch}{{\tilde{\mathcal L}}}
\renewcommand{\Mc}{{\mathcal M}}

\newcommand{\Tc}{{\mathcal T}}
\newcommand{\Uc}{{\mathcal U}}
\newcommand{\Vc}{{\mathcal V}} 
\newcommand{\Wc}{{\mathcal W}} 
\newcommand{\Xc}{{\mathcal X}}
\newcommand{\Yc}{{\mathcal Y}}
\newcommand{\Zc}{{\mathcal Z}}

\newcommand{\Cs}{\sc\mbox{C}\hspace{1.0pt}}

\newcommand{\Ds}{\sc\mbox{D}\hspace{1.0pt}}
\newcommand{\Ls}{\sc\mbox{L}\hspace{1.0pt}}

\newcommand{\Ts}{\sc\mbox{T}\hspace{1.0pt}}
\newcommand{\Ss}{\sc\mbox{S}\hspace{1.0pt}}
\newcommand{\Vs}{\sc\mbox{V}\hspace{1.0pt}}
\newcommand{\Us}{\sc\mbox{U}\hspace{1.0pt}}
\newcommand{\ah}{\tilde{\aF}}

\newcommand{\al}{\langle}
\newcommand{\ar}{\rangle}
\newcommand{\dist}{{\rm dist}\,}

\newcommand{\grad}{{\rm grad}\,}

\newcommand{\vol}{{\rm vol}\,}

\newcommand{\supp}{\hbox{{\rm supp}}\,}

\newcommand{\ND}{N \!\!\! D}
\newcommand{\DN}{D \!\!\! N}

\newcommand{\diam}{\operatorname{diam\,}}

\newcommand{\Sib}{\pmb{\pmb{\Si}}}
\newcommand{\Fib}{\pmb{\Phi}}

\newenvironment{dedication}{\itshape\center}{\par\medskip}

\newtheorem{thm}{Theorem}[section]
\newtheorem{cor}[thm]{Corollary}
\newtheorem{lem}[thm]{Lemma}

\newtheorem{proposition}[thm]{Proposition}

\theoremstyle{definition}
\newtheorem{defin}[thm]{Definition}

\newtheorem*{rem}{Remark}

\numberwithin{equation}{section}

\usepackage[hyperfootnotes=false,colorlinks=true,allcolors=blue]{hyperref}

\begin{document}


\title[Steklov eigenvalues Weyl asymptotics]{Weyl asymptotics for Poincar\'{e}-Steklov eigenvalues in a domain with  Lipschitz boundary}

\author{Grigori Rozenblum}
\address{Chalmers Univ. of Technology; The Euler International Mathematical Institute and St.Petersburg State University}
\email{$\mathrm{grigori@chalmers.se}$}
\begin{abstract}We justify the Weyl asymptotic formula for the eigenvalues of the Poincar\'e-Steklov spectral problem for a domain bounded by a  Lipschitz surface.
\end{abstract}
\thanks{The work of G.R. was performed at the Saint Petersburg Leonhard Euler International Mathematical Institute and supported by the Ministry of Science and Higher Education (Agreement No. 075--15--2022--287).}
\maketitle

\begin{dedication}
In the memory of Mikhail Birman and Mikhail Solomyak, my teachers.
\end{dedication}
\section{Introduction}\subsection{The Poincar\'e-Steklov problem.}
Let $\Om\subset \R^{d+1}$ be a connected bounded open set. We suppose that the boundary $\Si=\partial \Om $ is connected. The classical  \emph{Poincar\'e-Steklov (P-S) eigenvalue problem} consists in the study of the spectral properties for the Laplacian in $\Om$ with the spectral parameter entering in the boundary condition,
\begin{equation}\label{Problem}
  -\D u(x)=0, \, x\in \Om; \partial_\n u(x)-\la^{-1} u(x)=0, \, x\in\Si,
\end{equation}
where $\partial_\n$ is the derivative in the direction of the exterior normal $\n(x)$ at $x\in\Si.$ Further on, we assume that the boundary $\Si$ is Lipschitz and we equip $\Si$ with the surface measure induced by the Lebesgue measure in $\R^{d+1},$ which coincides with the $d$-dimensional Hausdorff measure.

This problem, often attributed to and named after V.A. Steklov,  \cite{Steklov}, was, in fact, first considered by H. Poincar\'{e} in 1896, see \cite{Poincare},  in relation to the analysis of tidal waves. Another early application, dealing with liquid waves in an open container was studied by D. Hilbert in \cite{Hilbert}.

The P-S spectral problem and its numerous generalizations found a lot of applications in physics and technology (see, e.g., \cite{AVKS} and the quite recent review  \cite{Developments}, for a, far from complete, bibliography on this topic.) Studies dealing with this problem continue up to now. The Steklov eigenvalue
problem appears in quite a few physical fields, such as
fluid mechanics, electromagnetism, elasticity, etc. It has applications for the study of various kinds of wave phenomena (see, e.g. \cite{Suslina21}), as well as in the seismology and tomography. Mathematically, this problem keeps being in the center of interest in spectral geometry and approximations.
MathSciNet shows more than 600 publications where various facets of Steklov type spectral problems are dealt with.

One of traditional topics in this field is the study of the  behavior of eigenvalues of the P-S problem.
 In the earliest paper by L. Sandgren on the asymptotics of these eigenvalues (see \cite{Sandgren}), inspired by Lars G{$\aa$}rding and {$\AA$}ke Plejel, with co-operation of Lars H\"ormander,  the author considered a moderately smooth domain $\Om$  in a Riemannian manifold.  The boundary of the domain was supposed be of the class $C^2$, but the coefficients of the main operator were set to belong to $C^1.$ Under these conditions, the asymptotic formula of H. Weyl type was obtained. Further on, the conditions imposed on the domain and coefficients were being gradually relaxed. Finally, in 2006, M.S. Agranovich \cite{Agr06} established the Weyl type formula for the case of a Lipschitz boundary, being, however, 'almost smooth' in the sense that this boundary should be at least  $C^1$ outside a closed set of zero surface measure; the symmetric second order elliptic operator $\Lc$ replacing the Laplacian in \eqref{Problem} is supposed to be of divergence form with continuous coefficients. Even earlier, probably, since \cite{AgrAmos}, M.S. Agranovich started popularizing the problem of extending the eigenvalue asymptotics results to domains with Lipschitz boundary, without additional assumptions. This question was morally supported by a flow of impressive results on various properties of elliptic boundary problems in Lipschitz domains by A.P.Calder\'{o}n, B.Dahlberg, C.Kenig, S.Hoffman, M.Taylor, M.Mitrea, E.Fabes, G.Verchota, M.S.Agranovich himself, and many others. It turned out that the setting of Lipschitz boundary problems and the typical  results are usually more or less the same as for a more smooth boundary, while the methods used in analysis  need to be essentially new -- and often quite complicated. This latter  circumstance is, partially, caused by the fact that certain important results in potential theory break down when passing from $C^1$ to Lipschitz surfaces.

 Meanwhile, the study of the P-S problem continued. We address the Readers to   \cite{Arendt}, \cite{ArendtMazz}, \cite{Behrndt},  \cite{Ouhabaz}, and references therein, where various  aspects of this problem were investigated. Important results were obtained even for boundaries considerably less regular than the Lipschitz ones, see, e.g., \cite{Taylor}.  Quite extensive became studies in the spectral geometry relating the geometric properties of the boundary to spectral properties. Especially rich were the results in the two-dimensional case. Here, it was known since long ago, see \cite{GR86}, that in the infinitely smooth case, the Steklov eigenvalues are, faster than any negative  power of their sequential number, close to those for the disk with the same perimeter. Quite a lot of further results in the two-dimensional case were  obtained since then. A description and huge bibliography can be found in \cite{Gorouard}, \cite{Developments}, and  \cite{SpecGeomBook}. We just mention here that in a series of papers, starting with \cite{Sloshing}, for the case of a piecewise smooth Lipschitz boundary, the authors succeeded in describing how the corners influence the deviation of the eigenvalues from their behavior in the smooth case. Meanwhile, in the  Lipschitz case, many important properties of the P-S (or the Neumann-to-Dirichlet, N-to-D,  $\ND$  operator) were established, see, e.g.,  \cite{Arendt}, \cite{ArendtMazz}, \cite{Behrndt}, \cite{Girou}; a number of important Steklov-type problems arising in hydrodynamics were considered recently in \cite{Suslina21}.

 Finally, after all these years, the essential progress in the eigenvalue asymptotics for a Lipschitz boundary was made  in \cite{LKP}, where, in the two-dimensional case, the Weyl eigenvalue asymptotics for the problem \eqref{Problem} was proved. Moreover,  in  \cite{LKP}, the boundary of the domain may be even somewhat more rough than Lipschitz. The method in  \cite{LKP} is based upon some deep results in the theory of conformal mappings, therefore, seemingly, it cannot be extended to higher dimensions. Thus, the eigenvalue asymptotics problem for Lipschitz boundary remained unresolved. In the recent  fundamental review \cite{Developments} and in \cite{Girou}, this problem is listed among unsolved and challenging.

 The author of the present paper, together with Grigory Tashchiyan, spent quite a lot of time attacking this problem, being initially inspired by M.S. Agranovich. Our hope was based upon our  result in \cite{RT1}, where, for integral operators of the type of single layer potential on a Lipschitz surface, a Weyl type eigenvalue asymptotic formula was proved. Since the P-S operator can be expressed in a simple way via the single layer and double layer potentials (see, e.g. \cite{Agr.UMN}, \cite{Agr}, and in a very general setting, \cite{Taylor}), it would be sufficient, for example, to know that the double layer potential is a compact operator. Unfortunately, for a Lipschitz surface, one should not expect this property of the double layer potential to hold; in fact, the opposite is known. Our efforts to circumvent this obstacle led to some partial results (see a mentioning in \cite{Suslina21}); we, however, did not consider these partial results deserving being published (a brief description of our efforts here, is given in Appendix, in hope that someone might be interested and more lucky -- if successful, such proof will, probably, be more esthetic than the present one).

A different approach is used in \emph{this} paper, not based directly upon the results in \cite{RT1} but rather on the \emph{idea } in \cite{RT1}, a quite natural one, of approximating the Lipschitz surface by  smooth surfaces and tracing how the corresponding compact Neumann-to-Dirichlet operators converge. In the study of solvability of elliptic boundary problems in Lipschitz domains, such approximation was, probably, first used by G.Verchota in \cite{Verchota}. In fact, a somewhat different realization of this approximation idea is implemented here. By a change of variables, the P-S problem in a Lipschitz domain for an elliptic operator in divergence form   with  coefficients continuous at the boundary,  is relocated to the problem for  another elliptic operator, this time in a smooth domain. Unfortunately, under such transformation, the coefficients of the operator \emph{cease to be continuous} even at the boundary, but stay  only bounded. Functions in $L_\infty$ cannot be approximated by smooth functions in $L_\infty$ norm. We, however, construct an approximation of the coefficients in a certain $L_p$-related norm, $p<\infty,$ and after a series of further transformations, we succeed in establishing a sufficiently strong convergence of operators describing the P-S spectrum, which enables us to perform the passage to the limit in eigenvalue asymptotic formulas.

  The reasoning in the paper  is, to a large extent, based upon the ideas the author absorbed many years ago while being a student of M.Sh. Birman and M.Z. Solomyak;  the paper is dedicated  to their memory.  The author  thanks Prof. Tatiana Suslina for  useful discussions and Prof. Jean Lagac\'{e} who acquainted him with an early version of the paper \cite{LKP}, which encouraged our efforts.  Permanent discussions with Prof. Grigory Tashchiyan for at least 15 years contributed a lot to a better understanding of the problem and stimulated the author, who is deeply grateful to Grigory for his involvement.

\subsection{Setting}
We are going to study the distribution of eigenvalues and singular numbers of various kinds of compact operators. By $n_{\pm}(\la, \Kb)$ we denote the  distribution function of positive (negative) eigenvalues of the \emph{self-adjoint} operator $\Kb$;  for an arbitrary compact operator $\Kb$,  $n(\la, \Kb)$ denotes the counting function of singular numbers of $\Kb.$ In the case when the operator is described by a quadratic form $\BFB[u]$ in the Hilbert space with norm defined by means of the quadratic form $\AF[u]$, we replace $\Kb$ by the ratio of these quadratic forms in the notation $n(\la,.)$, $n_{\pm}(\la,.)$.  For $\theta>0$, we denote by $\nb^{\sup}_\pm(\theta, \Kb)$ the quantity $\limsup_{\la\to 0}\la^{\theta}n_{\pm}(\la, \Kb).$ Without the subscript $\pm$, these notations concern the distribution of singular numbers of these operators. Without the superscript $\sup$,  these symbols denote the limits $\lim_{\la\to 0}\la^{\theta}n_{\pm}(\la, \Kb),$ provided these limits exist. Finally, if some operator or the ratio of quadratic form is defined in the displayed formula with tag (X.Y), the above notation  is used with this tag replacing the notation of the operator, e.g., $\nb^{\sup}(\theta, (X.Y))$, etc.

We consider  a bounded domain $\Om\subset\R^{d+1}$ with connected Lipschitz boundary $\Si$. (The connectedness condition can be easily removed, at the cost of certain notational complications.) This means that $\Sigma=\partial\Om$ can be covered by a finite collection of co-ordinate neighborhoods $\Uc_{\io}$ such that the portion of $\Si$ in $\Uc_\io$ can be, in a properly rotated  Euclidean co-ordinate system, $(x', x_{d+1})\in\Uc_\io$, represented by the equation $x_{d+1}=\psi^{(\io)}(x'), \, x'\in \Vc\subset\R^{d}$ with Lipschitz function $\psi^{(\io)}$ (sometimes such surfaces are called strongly Lipschitz in the literature, see, e.g., \cite{HMT}.) It is convenient to use a global co-ordinate system in a neighborhood of the boundary, see Sect.3.

In \cite{Sandgren}, \cite{Agr06}, and some other papers, a more general setting is considered, namely the spectral problem \eqref{Problem}, with the Laplacian in the equation replaced by some formally self-adjoint second order elliptic differential operator $\Lc$ in divergence form,  with the normal derivative $\partial_\n$ replaced by the derivative along the conormal  associated with $\Lc$ and with a weight function present on the right-hand side of the eigenvalue equation.  For us, it is also convenient to consider such more general, the \emph{weighted} $\ND$, Neumann-to-Dirichlet, problem. Namely, for a  uniformly positive definite (elliptic) real matrix-function  $\aF(x)=(a_{j,k}(x), \, j,k=1,\dots, d+1)$, $\aF\in L_{\infty}(\Om),$  we denote by $\Lc=\Lc_{\aF}$ the formal differential operator
\begin{equation}\label{L}
  \Lc\equiv\Lc_{\aF}=-\sum_{j,k}\partial_j a_{jk}(x)\partial_k.
\end{equation}
The weighted $\ND$ problem for $\Lc$ is the eigenvalue problem
\begin{equation}\label{LSteklov}
  \Lc u(x)=0, \, x\in\Om; \partial_{\aF}u(x)=\la^{-1} \ro(x)u(x), \, x\in\Si;\, \partial_{\aF}=\sum_{j,k}a_{j,k}\n_k\partial_{x_j},
\end{equation}
with a real function $\ro(x).$ The conormal  differential operator $\partial_{\aF}$ has symbol $\imath\langle \aF(x)\x,\n_x\rangle.$

A direct setting of the problem \eqref{LSteklov} requires an explicit description of the space of functions $u$ where the equation is considered. Under very weak conditions imposed here  on $\Om$ and $\aF,$ this task is very hard. Instead of this, we follow \cite{Sandgren}, \cite{Agr06} in considering the $\ND$ problem in the variational form, where all spaces under consideration admit explicit description, see Section 2. It is noted in \cite{Agr06} that such variational approach follows the way of reasoning of H.Poincar\'{e}, V.A.Steklov and other researchers of that time, for whom it is the variational setting of the eigenvalue problem, based upon considering energy balance, was the primary point of the eigenvalue analysis, while the differential Euler-Lagrange equation appeared as a secondary object.

\subsection{Asymptotic formulas. The main result}
The asymptotic eigenvalue formula for the problem \eqref{L} can be written in the following way (see, e.g., \cite{Agr06}). In a fixed orthogonal  co-ordinate system, we associate with the matrix $\aF(x)$ its sesquilinear form $\pmb{\aF}_x$: $\pmb{\aF}_x(\x)=\sum_{j,k} a_{jk}(x)\x_j\x_k;$ let $\pmb{\aF}_x(\x,\y)$ be the corresponding bilinear form.
For $x\in\Si,$ $\x'\in \mathrm{T}^*_{x}\Si$ and the normal vector $\n$ at $x,$ we denote by $\beta(x;\x')$ the positive square root of $\pmb{\aF}_x(\n)\pmb{\aF}_x(\x')-\pmb{\aF}_x(\x',\n)^2$,

\begin{equation}\label{beta}
  \beta(x,\x')=(\pmb{\aF}_x(\n)\pmb{\aF}_x(\x')-\pmb{\aF}_x(\x',\n)^2)^{\frac12};
\end{equation}
this is a function positively homogeneous in $\x'$ of order $1$.
 Then the density $\a_{\pm}(x)$ is defined as
\begin{equation}\label{alpha}
  \a_{\pm}(x)=\vol\{\x'\in \mathrm{T}^*_{x}\Si:\beta(x;\x')<\ro_{\pm}(x) \},
\end{equation}
and the asymptotic formulas to be proved are
\begin{equation}\label{formula}
  \nb_{\pm}(d,\ref{LSteklov})=(2\pi)^{-d}\int_{\Si}\a_{\pm}(x)d\m_{\Si}(x),
\end{equation}
where $\m_{\Si}$ is the surface measure on $\Si$ generated by the embedding of $\Si$ into $\R^{d+1}.$

For a \emph{continuous} matrix-function $\aF,$ such definition of the eigenvalue problem \eqref{LSteklov} and the understanding of the formula \eqref{formula} do not cause confusion, even if the matrix $\aF$ is continuous only at the boundary $\Si.$ In our case,  the matrix $\aF(x)$ may be discontinuous. However, under the conditions imposed below, the limit values of $a_{jk}$ at the boundary are well defined almost everywhere in $\Si,$ and therefore these formulas \eqref{formula} make sense.

In \cite{Agr06}, the asymptotic formula \eqref{formula} was justified under the following conditions, see \textbf{Theorem} on P.242.
\begin{thm}\label{Th.Agr} Let $\Om\subset \R^{d+1}$ be a bounded domain; suppose that $\Si=\partial\Om$ is a surface  of class $C^1$ and $\ro\in L_\infty(\Si).$  Suppose finally that  the coefficients $a_{jk}$ are continuous in $\overline{\Om}$. Then  the asymptotics  \eqref{formula} holds. Moreover, the same conclusion is correct if $\Si$ is Lipschitz but belongs to $C^1$ in a neighborhood  of the support of $\ro,$  outside a closed set of zero surface measure. The function $\ro$ may belong to $L_d(\Si)$ for $d>1$ and any $L_q(\Si)$ with  $q>1,$ for $d=1.$
\end{thm}


The \emph{elementary} formulation of the main result of our present paper  is the following.

\begin{thm}\label{RTthm} Let $\Om\subset\R^{d+1}$ be a bounded domain with Lipschitz boundary $\Si.$ Suppose that the operator $\Lc_\aF$ is uniformly elliptic in $\Om$, the coefficients $a_{jk}$ are continuous in $\overline{\Om}$. Suppose that  $\ro\in L_d(\Si)$ for $d>1$ and $\ro$ belongs to the Orlicz class $L\log L(\Si)$ for $d=1.$ Then the asymptotic formula \eqref{formula} is valid.
\end{thm}
The continuity condition in Theorem \ref{RTthm} can be considerably relaxed. We say that a function $a(x)\in L_\infty(\Om)$ is \emph{continuous at the boundary} if there exists a continuous function $a_b\in C(\Si),$ $\Si=\partial\Om$ such that for any $x'\in\Si$,
\begin{equation}\label{Def.Cont.bdry}
  \lim_{\de\to 0}\|a(.)-a_b(x')\|_{L_\infty(B(x',\de)\cap\Om)}=0,
\end{equation}
where $B(x',\de)$ is the ball in $\R^{d+1}$ with radius $\de$, centered at $x'.$
Since the functions in $L_\infty(\Om)$ are defined up to a set of measure $0$, we may identify in the above definition the continuous function $a_b$ with the restriction of $a$ to $\Si$.

Having  this definition, the  formulation of the more general result that  we are going to prove here, is   as follows.
\begin{thm}\label{Thm.General}Let $\Om\subset\R^{d+1}$ be a bounded domain with Lipschitz boundary $\Si.$ Suppose that the operator $\Lc_\aF$ is uniformly elliptic in $\Om$, the coefficients $a_{jk}$ are bounded and measurable and they are continuous at the boundary in the sense of \eqref{Def.Cont.bdry}. Then the asymptotic formula \eqref{formula} is valid.
\end{thm}
\begin{rem}The conditions imposed on the coefficients of the operator $\Lc$ can be further relaxed, allowing them to be unbounded in a certain sense and admitting certain degeneration of ellipticity. We do not pursue this line here, in order to keep the elementary level of the presentation.
\end{rem}
We describe here briefly the structure of the proof. After certain transformations, the weighted $\ND$ problem in the domain with Lipschitz boundary is relocated to a similar  problem  for a  uniformly elliptic operator $\Lc$ but with bounded coefficients, possibly, discontinuous  at the boundary,  in a domain with smooth boundary. This operator is approximated by operators $\Lch$ with smooth coefficients which converge to the coefficients of $\Lc$  in a certain $L_p$ - based norm, with some $p<\infty$. Then it is proved that the compact operators describing the P-S spectrum for $\Lch$ converge so strongly that it becomes possible to pass to the limit in the formulas for coefficients $\nb(d,.)$ in the Steklov eigenvalue asymptotics for $\Lch.$

\section{The variational representation}\label{Sect.Var}
The setting of the eigenvalue problem under  very weak regularity conditions imposed on the coefficients and the domain, is variational. As usual, for more regular data, the variational setting is equivalent to the classical, 'strong' one.
\subsection{The weighted $\ND$ operator}
Since the earliest paper on the topic by Lenart Sandgren \cite{Sandgren}  in 1955, the  most commonly used method of the study of Steklov eigenvalues is the variational one. Its main advantage is its robustness with respect to various perturbations of the problem. We too will  need a version of the variational setting of the Steklov eigenvalue problem (see, however, Appendix where we discuss a different approach).

First of all, since we are going to use the instruments of the perturbation theory for compact operators, we will need to pass from the unbounded Dirichlet-to-Neumann operators $\DN$ (which are often considered) to the N-to-D ones, $\ND,$ since the latter operators are compact (this corresponds to the placement of the spectral parameter in (1.3)). A minor inconvenience here is that the Neumann problem for the Laplacian (as well as for other second order  elliptic operators without zero order terms) has an eigenfunction with zero eigenvalue, so the Neumann operator is not invertible. There are several standard ways to deal with this  circumstance. Say, often, the Neumann operator is restricted to the codimension 1 subspace of functions orthogonal to constants. In our considerations, it is more convenient  to avoid hitting the zero eigenvalue by means of adding a nontrivial nonnegative  function $\vb_0(x)$  to the operator $\Lc$, namely to consider the operator $\Lc+\vb_0(x)$. Such perturbation is relatively compact with respect to  self-adjoint realizations of $\Lc$, and, by the usual perturbation reasoning, see, e.g., Lemma 1.3 in  \cite{BS}, it does not influence the asymptotics of the spectrum. One might have set $\vb_0\equiv 1,$ but this does not simplify the matter. For  economy of symbols, we use the notation $\Lc$ for the operator with $\vb_0$ already absorbed, further on. At a proper moment, we will use the freedom in choosing $\vb_0.$

We denote by $\Tb\equiv\Tb_{\aF}$ the operator $\Lc$ with \emph{Neumann} boundary conditions on $\Si$. Namely, it corresponds to the quadratic form
\begin{equation}\label{form ab}
\ab_0[u]=\int_\Omega( \ap [ u](x) +\vb_0(x)|u(x)|^2) dx, \, \ap [u](x):=\sum_{j,k} a_{j,k}(x)\partial_j u(x) \partial_k \bar{u}(x),
\end{equation}
with domain $u\in H^1(\Omega); $ this form is  equivalent to the square of the standard norm $\|u\|_{H^1(\Om)}$ in the Sobolev space $H^1(\Omega).$ Being considered in $L_2(\Om)$, this quadratic form defines the positive self-adjoint operator $\Tb$ satisfying
 \begin{equation*}
  \ab_0[u]=(\Tb^{\frac12}u,\Tb^{\frac12}u), \, u\in H^1(\Om)=\Dc(\Tb^{\frac12}).
 \end{equation*}
(This fact, which we will use systematically, that the domain of a closed positive  quadratic form coincides with the domain of the square root of   the operator defined by this quadratic form, is a standard property in the abstract spectral theory, see, e.g., \cite{Kato}, Theorem 2.23 in Sect. VI.  For quadratic forms which are only sectorial (the case we do not need), this equality constitutes the famous \emph{Kato's square root problem} which has been the topic of deep extensive studies since 1980-s.)

{Further on,  we will sometimes use the notion of the spectrum of the ratio of quadratic forms in the Hilbert space. Actually, the spectrum of the ratio $\frac{\BFB[f]}{\AF[f]}$ with $f\in\XF$ is a shorthand for the spectrum of the operator defined by the quadratic form $\BFB[f]$ in the Hilbert space $\XF$ equipped with the norm defined by the quadratic form $\AF[f]$. }


For a real function $\ro_0(x), x\in \Si$ we consider the quadratic form
\begin{equation}\label{numerator}
  \rob_0[f]=\int_{\Si}\ro_0(x)|f(x)|^2 d\m_{\Si}(x).
\end{equation}
Let $\ro_0$ belong to $ L_d(\Si).$
By Kondrashov's trace theorem, for $d>1,$ the trace mapping $u\mapsto u_{|_\Si},$ defined initially on continuous functions  $u\in H^1(\Om)$, extends by continuity to the bounded mapping $\g: u\mapsto u_{|_\Si},$ $\g: H^1(\Om)\to L_q(\Si),$ $q=\frac{2d}{d-1}. $ Therefore, by the H\"older inequality, if $\ro_0\in L_d(\Si)$, with $d^{-1}+2 q^{-1}=1$, the quadratic form \eqref{numerator} satisfies
\begin{equation}\label{traceform}
\rob_0[\g u]\le C \|\ro_0\|_{L_{d}(\Si)}\ab_0[u], \, u\in H^1(\Om).
\end{equation}
For $d=1,$ a somewhat more complicated reasoning (see, e.g., \cite{AdF}) gives
\begin{equation}\label{traceform1}
  \rob_0[\g u]\le C \|\ro_0\|_{L\log L(\Si)}\ab_0[u], \, u\in H^1(\Om),
\end{equation}
where ${L\log L(\Si)}$ is the Orlicz space consisting of functions $u$ on $\Si$ for which $|u|(1+|\log|u||)\in L_1(\Si),$ with the Luxemburg norm.

The following variational representation is known since \cite{Sandgren}, see e.g.,  \cite{Agr06}, \cite{Taylor},  and sources cited there: the Steklov spectrum is described by the ratio of quadratic forms $\frac{\rob_0[\g u]}{\ab_0[u]}$ considered on the space of $L_2-$solutions of the equation $\Lc u=0.$ Since $\rob_0[\g u]=0$ for $u\in \overset{\circ}{H}{}^1(\Om),$ one can, as this is usually done, drop the condition $\Lc u=0$ in the variational setting and consider the spectrum of the ratio
\begin{equation}\label{ratio}
  \frac{\rob_0[\g u]}{\ab_0[u]}, \, u\in H^1(\Om).
\end{equation}
\subsection{A freedom in the choice of the weight function $\ro_0$.}
Estimates \eqref{traceform}, \eqref{traceform1} mean that the quadratic form \eqref{numerator} is bounded in $H^1(\Om)$,  the equality in \eqref{numerator} can be therefore extended by continuity to all traces of functions in $H^1(\Om)$ and therefore the form $\rob_0$ defines a bounded operator $\Sbb=\Sbb[{\ro_0}]$ in $H^1(\Om)$, with norm estimated by $\|\ro_0\|_{L_d(\Si)}, $ resp., $\|\ro_0\|_{L\log L(\Si)}.$ If the function $\ro_0$ is real-valued, this operator is self-adjoint.
In fact, this operator is compact and its eigenvalues satisfy a power order estimate.

\begin{thm}\label{Th.estimate} Let $\Om\subset\R^{d+1},$ $d>1,$ be a bounded  domain with Lipschitz boundary $\Si$ and let $\ro_0\in L_d(\Si)$. Then for the singular numbers of $\Sbb=\Sbb{[\ro_0]}$ the estimate holds
\begin{equation}\label{Est}
  n(\la,\Sbb{[\ro_0]})\le C(\Om,\aF)\la^{-d}\|\ro_0\|_{L_d(\Si)}^{d}.
\end{equation}
If the function $\ro_0$ is real-valued, estimates of the form \eqref{Est} hold, separately, for positive and negative eigenvalues,
\begin{equation}\label{Est.pm}
  n_{\pm}(\la,\Sbb{[\ro_0]})\le C(\Om, \aF)\la^{-d}\|(\ro_0)_{\pm}\|_{L_d(\Si)}^{d}.
\end{equation}
For $d=1$ the above results hold, with $L_d(\Si)$ replaced by $L\log L(\Si).$
\end{thm}
\begin{proof} For $d>1,$ we  use the estimate for the eigenvalues of Birman-Schwinger type operators with singular measures, obtained recently in \cite{RTsing}. We consider some bounded open set $\Om'$ containing $\Om$ strictly inside and apply Theorem 3.3 in \cite{RTsing}, for the particular choice of $\Nb=d+1, l=1, s=d, d>1, $ $\m$ being the surface measure on the surface $\Si$. This theorem states the following. Let $\Om'\subset\R^{\Nb}$ be a bounded open set, $\m$ be a compactly supported measure in $\Om'$ satisfying, for  $X\in\Si,$ $r\le \diam(\Si)$, $\Si=\supp\m$,
\begin{equation}\label{GTest}
  \m(B(X,r))\le \Ac r^{s},\, s>\Nb-2(=d-1), \ro\in L_{\theta,\m}, \, \theta =\frac{s}{s+2-\Nb},
\end{equation}
where $B(X,r)$  is the ball with radius $r$ centered at $X.$  Then, for the particular case  $s=d,$ $\theta=d,$  for the operator $\Sbb=\Sbb[\ro_0,\Om']$
in $\overset{\circ}{H}{}^1(\Om')$ defined by the quadratic form $\rob_0$ in \eqref{numerator}, the estimates \eqref{Est}, \eqref{Est.pm} hold.

We use this result in the following way. First, note that for a compact Lipschitz surface $\Si$ of dimension $d$ with $\m$ being the  $d$-dimensional Hausdorff measure on $\Si$, the condition \eqref{GTest} is satisfied with $s=d$. Next, by the Calder\'{o}n-Stein extension theorem, (see, e.g., Theorem 5.24 in \cite{AdF}), there exists a bounded extension operator $\Ec:H^1(\Om)\to \overset{\circ}{H}{}^1(\Om'),$
\begin{equation*}
\|\Ec v\|_{\overset{\circ}{H}{}^1(\Om')}\le C(\Om,\Om')\|v\|_{H^1(\Om)}.
\end{equation*}
Therefore,
if on some subspace $\Ls\subset H^1(\Om)$ of dimension $n,$ we have
\begin{equation*}
  \frac{\rob_0[\g v]}{\|v\|^2_{H^1(\Om)}}\ge  \la, \, v\in\Ls,
\end{equation*}
it follows that
\begin{equation*}
  \frac{\rob_0[\g \Ec v]}{\|\Ec v\|^2_{\overset{\circ}{H}{}^1(\Om')}}=\frac{\rob_0[\g  v]}{\|\Ec v\|^2_{\overset{\circ}{H}{}^1(\Om')}}\ge  C(\Om,\Om')^{-2}\la .
\end{equation*}
By the variational principle, this implies the inequality for the counting functions of eigenvalues of operators,
\begin{equation*}
  n_+(\la,\Sbb[\ro_0,\Om])\le n_+(C(\Om,\Om')^{-2}\la,\Sbb[\ro_0,\Om']),
\end{equation*}
which gives \eqref{Est.pm} with '+' sign. Other estimates follow in a similar way.\\
For $d=1,$ the reasoning is the same, just we use the Orlicz estimate in \cite{RSh} instead of the eigenvalue estimate in \cite{RTsing}.
\end{proof}

\begin{rem}In fact, some stronger results which, however, are not needed for our present topic, hold. Actually, the extension theorem for the class $H^1(\Om)$, which we used,  is valid under somewhat weaker restrictions than the Lipschitz property.  Namely, in \cite{Jones}, the class of $(\e-\de)$ - domains has been introduced (see the definition in \cite{Jones}, p.73) called also  \emph{locally uniform domains} in the literature. This class contains all Lipschitz domains, but some other, less regular ones, as well (it is mentioned in \cite{Jones} that an $(\e-\de)$-domain in $\R^{d+1}$ may have a fractal boundary of any Haussdorff dimension in $[d,d+1)$, particular examples being the Von Koch snowflakes. On the other hand, an exterior cusp  prevents the domain from being locally uniform).  For an $(\e-\de)$ domain, a bounded extension operator $\Ec$ for Sobolev spaces exists. Therefore, if the $d$-dimensional Hausdorff measure on $\partial\Om$ satisfies \eqref{GTest} with $s=d$ and $\Om$ is an $(\e-\de)$-domain then for the operator $\Sbb [\ro_0,\Om]$ with a function $\ro_0\in L_d(\m)$, which can be considered as a natural generalization of the $\ND$ operator, the estimate \eqref{Est}, \eqref{Est.pm} is valid.
\end{rem}

The operator $\Sbb\equiv\Sbb{[\ro_0]}$ has the same eigenvalues as the
 operator of the P-S problem in the smooth case. In the course of our reasoning, we will introduce some more operators with the same spectrum.

The eigenvalue estimates in Theorem \ref{Th.estimate}, with coefficient in front of the power term depending on the integral norm of the functional parameter, enable one to establish the property of the passage to limit in asymptotic eigenvalues formulas. This way of reasoning, starting in the classical studies by M.Sh. Birman and M.Z. Solomyak in 1960-s -1970-s (see, e.g., \cite{BS}, Lemma 1.5), is now the standard tool in establishing the eigenvalue asymptotics for singular problems by means of  approximating by more regular ones. We reproduce here this extraordinary lemma, which we use at least on three occasions in the course of the paper.

\begin{lem}\label{Lem.5.1}Let $\Kb$ be a compact operator, $\theta>0,$ and for sufficiently small $\e$ one can split $\Kb$ into the sum, $\Kb=\Kb_\e+\Kb_\e'$ so that for $\Kb_\e$  the singular numbers asymptotic formula holds,
\begin{equation}\label{5.1.1}
  \nb(\theta,\Kb_\e)=A_\e,
\end{equation}
 and the operators $\Kb_\e'$  are asymptotically small in the sense
 \begin{equation}\label{5.1.2}
   \nb^{\sup}(\theta, \Kb_\e')\to 0,
 \end{equation}
 as $\e\to 0$. Then the limit $\lim_{\e\to 0}A_\e=A$  exists and the asymptotics holds,
 \begin{equation}\label{5.1.3}
  \nb(\theta,\Kb)=A.
 \end{equation}
 An analogous statement is valid for the positive/negative eigenvalues of a \emph{self-adjoint} operator $\Kb$. One should replace in the above formulation $\nb$ by $\nb_{\pm}$ in \eqref{5.1.1} and \eqref{5.1.3} (but not in \eqref{5.1.2}).
\end{lem}

 Among most recent  examples of applying this way of reasoning, one can cite  \cite{RTsing}, Sect.6 and \cite{LKP}. We formulate here the particular statement concerning the eigenvalue asymptotics for the operator $\Sbb[{\ro_0}];$ it follows from estimates \eqref{Est}, \eqref{Est.pm}
 \begin{cor}\label{CorLimAs}
   Suppose that for the weights $\ro_0$ in some set $\Yc\subset L_{d}(\Si)$, $d>1,$ the asymptotic formula \eqref{formula} is established. Then this formula is valid for all $\ro_0\in\overline{\Yc},$ the closure of $\Yc$ in the norm of $L_{d}(\Si).$ The same statement is valid for the dimension $d=1,$ with $L_d(\Si)$ replaced by the Orlicz space $L\log L(\Si).$
 \end{cor}
 In particular, this means that we can restrict ourselves to considering as $\Yc$, the set of  continuous, and, later, smooth functions $\ro_0$, which are dense in the corresponding space with integral norm.

 We follow now \cite{Sandgren}, \cite{Suslina99}, \cite{VuSolIzv}, \cite{Agr}, \cite{RT1} and other papers where the variational method was used for the study of spectral problems containing a weight, including the  P-S type spectral problems. Namely, once   order sharp eigenvalue \emph{estimates}  involving an integral norm of the weight function are obtained, it is possible to restrict the further study  of eigenvalue \emph{asymptotics }to considering \emph{sign-definite} weight functions only. We do not need to repeat this standard  reasoning in detail, but just describe briefly  the common scheme. Starting with a non-sign-definite weight function $\ro_0$, we first approximate it in a proper integral  norm by a function $\ro$ which has 'separated' positive and negative parts, $\dist(\supp \ro_+, \supp \ro_-)>0.$ After this, the study of the distribution of the positive, resp., negative spectrum of the problem is performed by cutting-away the part of the weight with wrong sign.

 Consequently, we restrict ourselves to non-negative weight functions $\ro_0$ further on.
 \subsection{The asymptotic coefficient}
 The next preparatory step consists in obtaining a convenient expression for the coefficient in the asymptotic formula for the eigenvalues of the P-S problem so that it is possible to pass to the limit in this expression as the coefficients of the operator converge in a proper sense.

 The asymptotic coefficient, known  to be correct in the smooth case and aimed for in the Lipschitz case,
 is given by \eqref{formula}, \eqref{beta}, \eqref{alpha}.

 For our further needs, it will be more convenient to use a somewhat different representation of the coefficient in \eqref{formula}.
Namely, for a given $x\in\Si,$ we can represent $\ap_x(\x)$, $\ap_x(\x,\y)$ as matrix products
 \begin{equation*}
  \ap_x(\x)=\x^* \aF(x)\x, \,\, \ap_x(\x,\y)=\x^*\aF(x)\y=\y^*\aF(x)\x,
 \end{equation*}
(where $\x^*$ denotes the row-vector, matrix-adjoint to the column-vector $\x$.) In these notations, we have
\begin{gather}\label{coeff.1}
                \beta(x,\x')^2=\n^*\aF(x)\n {\x'}^*\aF(x)\x'- {\x'}^*\aF(x)\n\n^*\aF(x)\x'             \\\nonumber
      = \n^*\aF(x)\n\langle\aF(x)\x',\x'\rangle -     \langle [(\n^*\aF(x))^* (\n^*\aF(x))]\x',\x'\rangle :=\langle \Theta(x)\x',\x'\rangle
                           \end{gather}
where $\Theta(x)=(\n^*\aF(x)\n) \aF(x)-\aF(x)\n\n^*\aF(x)$ (note that the last matrix here has rank 1, together with  $\n\n^*$.).  Now let $\Pi_x$ be the embedding  of $\mathrm{T}^*_{x}\Si$ into $\R^{d+1},$ and  $\Pi_x^*$ be the adjoint operator, the  projection  of $\R^{d+1}$ to  $\mathrm{T}^*_{x}\Si.$
Then  \eqref{coeff.1} can be written as
  \begin{equation}\label{coeff.3}
  \beta(x,\x',\n)^2=\langle{\Pi_x^*\Theta_x \Pi_x}\x',\x'\rangle,
  \end{equation}
  therefore,
\begin{equation}\label{coeff.2}
  \a_{\pm}(x)=\vol\{\x'\in\R^d: 0<\langle \Theta'_x)\x',\x'\rangle^{\frac12}< (\ro_0)_{\pm}(x)\}=\om_{d}(\ro_0)_{\pm}(x)^{d}\det(\Theta'_x)^{-1/2},
\end{equation}
with $d\times d$ matrix $\Theta'_x=\Pi_x^*\Theta_x\Pi_x$ The meaning of the  matrix $\Theta'_x$ can be understood in the following way. We represent the $(d+1)\times (d+1)$ matrix $\Theta$ in an orthogonal frame where one of basis vectors is the normal vector $\n_x.$ Next, we strike out in this matrix the row and the column corresponding to $\n_x$. What remains is just our matrix $\Theta'_x$.
Using \eqref{coeff.2}, we can  establish certain estimates for the quantity $\nb_{\pm}(d,\ref{LSteklov}),$ uniform in a class of elliptic operators.

It follows from the Cauchy-Schwartz inequality that
\begin{equation}\label{unif.1}
  \pmb{\aF}_x(\n)\pmb{\aF}_x(\x')-\pmb{\aF}_x(\x',\n)^2=\langle\aF(x) \n,\n\rangle\langle\aF(x) \x',\x'\rangle-\langle\aF(x) \x',\n\rangle^2\ge 0.
\end{equation}
For a positive matrix $\aF,$ the exact equality in the Cauchy-Schwartz inequality \eqref{unif.1} is possible only if the vectors $\x'$ and $\n$ are parallel. The latter may never happen, therefore the inequality in \eqref{unif.1} is strict; by homogeneity we  obtain
\begin{equation}\label{unif.2}
  \pmb{\aF}_x(\n)\pmb{\aF}_x(\x')-\pmb{\aF}_x(\x',\n)^2\ge \kb(\x'/|\x'|)\pmb{\aF}_x(\x'),
\end{equation}
with some $\kb(\x'/|\x'|)>0.$

The sharp constant $\kb(\x'/|\x'|)$ in \eqref{unif.2} depends continuously on $\x'\ne 0;$ by the compactness of the unit sphere, we have
\begin{equation}\label{unif.3}
  \be_x(\x')^2=\pmb{\aF}_x(\n)\pmb{\aF}_x(\x')-\pmb{\aF}_x(\x',\n)^2\ge \kb(\x'/|\x'|)\pmb{\aF}_x(\x')\ge C\kF(x)|\x'|^2,
\end{equation}
where $\kF(x)$ is the ellipticity constant of the matrix $\aF(x)$ at the point $x\in\Si.$ We suppose that the matrix $\aF(x)$ is  uniformly elliptic, $\aF(x)\ge \kb_0>0.$ Then \eqref{unif.3}  implies
\begin{equation}\label{unif.4}
 \a_{\pm}(x)\le C\kb_0^{-d}\om_d\ro_{\pm}(x)^d.
\end{equation}
In the equivalent representation, the estimate for the matrix  $\Theta'_x$  also follows, $\det(\Theta'_x)\ge C\kb_0^{-d}.$

The last estimates enable us to establish the  following convergence property.
\begin{lem}\label{lem.convergence}
Let $\aF_s(x)$, $s\in[0,1)$ be a family  of Hermitian matrix functions on $\Si$, such that they satisfy ellipticity estimates $\aF_s(x)\ge \kb_0$  uniformly in $x, s,$   and are also uniformly  bounded, $|\aF_s(x)|\le \kb_1$. Suppose that the matrices $\aF_s$ converges as $s\to 0$ to $\aF_0$ in $L_q(\Si),$ $q<\infty.$ Denote by $\varrho_s(x)$ the function
 \begin{equation*}
  \varrho_s(x)=\det(\Theta_s'(x))^{-\frac12},
 \end{equation*}
where $\Theta_{s}'$ is the matrix $\Theta'$ in \eqref{coeff.1}, corresponding to the coefficient matrix $\aF_s$. Then
   \begin{equation*}
   \varrho_s(x)\to\varrho_0(x)
   \end{equation*}
 in $L_q(\Si).$
\end{lem}
\begin{proof}
  The statement follows easily from the fact that $\varrho_s(x)^2$ is a rational function of the entries of the matrix $\aF_s(x)$ and these matrices are  separated from zero.
\end{proof}
The following  result shows that the expression for  the coefficient in the asymptotics of P-S eigenvalues endures the convergence of the coefficient matrix in the integral metric.
\begin{lem}\label{limitformula}
  Let $\aF_s(x)$ be a family of matrices, as in Lemma \ref{lem.convergence} and $\ro_0$ be a bounded function on $\Si.$
  Then  $\int \varrho_{s,\pm}(x)\ro_0(x)d\Si\to\int \varrho_{\pm}(x)\ro_0(x)d\Si.$
\end{lem}
The statement follows by the passage to limit, using Lemma \ref{lem.convergence}.
\subsection{D-to-N and N-to-D operators as pseudodifferential ones}
It is well known that for the Laplacian, the $\DN$ operator equals, up to lower order terms, $(-\Delta_{\Si})^{\frac12}$, for a smooth boundary; an interesting  discussion can be found in \cite{Girou}. For general elliptic operators, in \cite{Agr.UMN}, \cite{AgrTMMO},  M.S.Agranovich presented an explanation of the coefficient $\beta(x,\x')$ for the case of a smooth (or 'almost smooth') surface $\Si.$  We discuss now the latter formula. Suppose first that the coefficients matrix $\aF(x)$ does not depend on $x\in\Om.$ Consider the fundamental solution $R(x-y)$ for the operator $\Lc,$ so that $\Lc_x R(x-y)=\delta(x-y).$ With this fundamental solution the classical potential operators are associated, namely, the single layer potential operator $\Ss$,
\begin{equation*}
f\mapsto  \Ss_f(x)=\int_{\Si}R(x-y)f(y)d\mu(y),\, x\in\Si,
\end{equation*}
and the 'direct value'  of the conormal derivative of the single layer potential (often called the Neumann-Poincar\'{e} operator)
\begin{equation*}
\Ds': \,f\mapsto \Ds'_f(x)=\int_{\Si}\partial_{\aF(x)}R(x-y)f(y)d\mu(y),\, x\in\Si,
\end{equation*}
the latter integral understood in the principal value sense.
The adjoint to $\Ds'_f$ is the integral operator $\Ds,$ called the direct value of the double layer potential,
 \begin{equation*}
 \Ds:\,f\mapsto \Ds_f(x)=\int_{\Si}\partial_{\aF(y)}R(x-y)f(y)d\mu(y),\, x\in\Si.
 \end{equation*}

Then the $\ND$ operator is, up to lower order terms, the composition
\begin{equation}\label{D2N}
  \ND =(\frac12+\Ds')^{-1}\Ss= \Ss(\frac12+\Ds)^{-1}.
\end{equation}
For a smooth boundary $\Si$, the operators $\Ss$ and $\Ds$ are order $-1$ pseudodifferential operators on $\Si.$ Since $\Ss$ is a restriction of the order $-2$ pseudodifferential operator to the boundary, its symbol is calculated according to the usual rules of pseudodifferential calculus,
\begin{equation}\label{D2N rest}
\ssF(x,\x')=\frac{1}{2\pi}\int_{-\infty}^\infty[\aF(\x',\x_{d+1})]^{-1}d\x_{d+1},\, x\in\Si,
\end{equation}
in co-ordinates where $x_{d+1}$ is directed along the normal vector  $\n_x$ to the tangent plane at the point $x\in\Si$ and the corresponding co-ordinates $\x',\x_{d+1},$  with $\x'$ in the cotangent plane and $\x_{d+1}$ directed along $\n_x$,
see, e.g., \cite{Agr.UMN},  (4.28). Calculations by $\eqref{D2N rest}$ give the expression for  $\be(x,\x')=\ssF(x,\x')^{-1}$ in \eqref{beta} as the principal symbol of the operator $\DN=\ND^{-1}.$
Similar considerations cover the case of variable smooth coefficients, where the fundamental solution takes the form $R(x,y),$
$\Lc_xR(x,y)=\delta(x-y)$ and in \eqref{D2N rest} one should change  $\aF(\x',\x_{d+1})^{-1}$ to $\aF(x',0,\x',\x_{d+1})^{-1},$ see \cite{Agr}, Sect.3.4.
For a non-smooth, Lipschitz surface, the function $\beta(x,\x')$ still can be calculated by the above formula, but it is not a symbol of a nice pseudodifferential operator any more.

The principal symbol of the double layer operator $\Ds$ can be also calculated, see, e.g., \cite{Taylor2}, but the we do not need an explicit expression here.

 \section{Some geometry considerations}\label{Sect.geom}  In this Section we introduce a global co-ordinate system in an interior collar neighborhood of the boundary $\Si=\partial\Om.$ After the passage to this co-ordinate system and a special co-ordinates change, the Poincar\'e-Steklov problem in a Lipschitz domain, for the equation with  coefficients, \emph{continuous at the boundary}, reduces to  an equivalent problem in a domain with smooth boundary, but for an elliptic operator with \emph{discontinuous}  bounded measurable coefficients, fortunately,   still with a certain very weak continuity property.

 \subsection{A smooth  surface approximating the boundary}
 By the Rademacher theorem, a Lipschitz function  has derivatives almost everywhere and these derivatives are bounded. Therefore, the Lipschitz  surface $\Si\subset\R^{d+1}$ has a tangent plane almost everywhere. This fact enables one to describe explicitly the surface measure on $\Si$ generated by the Lebesgue measure on $\R^{d+1}.$ Namely, on the local Lipschitz graph $x_{d+1}=\psi(x'), \, x'\in\Vs\subset\R^{d},$ the surface measure is given by
 \begin{equation}\label{Lip1}
   \m_{\Si}(E)=\int_{\pi E}\sqrt{1+ |\nabla\psi(x')|^2}dx',
 \end{equation}
where $\pi E$  is the projection of the set $E\subset \Si$ to $\R^{d}.$  It is important to note that for a Lipschitz function $\psi,$ the integrand  in \eqref{Lip1} is a bounded function, an algebraic expression of the partial derivatives of $\psi.$

A detailed study of geometric properties of Lipschitz domains is presented in \cite{HMT}. In particular, for a domain $\Om\subset \R^{d+1}$ with Lipschitz boundary, the distributional gradient $\grad(\hi_\Om)$ of the characteristic function of $\Om$ is a vector measure which is absolutely continuous with respect to the measure $\m_{\Si}$ in \eqref{Lip1},

\begin{equation}\label{Lip2}
  \grad(\hi_\Om)=-{\n}(x)\m_{\Si}
\end{equation}
with the vector-valued density ${\n}$, $|{\n}(x)|=1$ coinciding with the normal at $x$ almost everywhere on $\Si.$
(Such vector measure can  be well defined for a class of domains even  somewhat less regular as well, see \cite{HMT}.)

Starting from \cite{Sandgren}, and then in   \cite{Agr06}, \cite{Suslina99}, \cite{Agr}, etc., a localization was used, reducing the spectral analysis for operators related with the P-S problem in $\Om$ to the ones in Lipschitz cylinders. In our  study, it is more convenient to define the surface $\Si$ and perform related constructions globally.

To do this, we will need the classical result on the existence of a smooth vector field transversal to the Lipschitz surface. Namely, there exists a $C^{\infty}$-vector field $\pmb{\G}(x),$ $x\in \R^{d+1}$ such that at the points $x\in\Si=\partial \Om,$ those points  where the tangent plane
exists, $|\pmb{\G}(x)|=1$ and $\pmb{\G}(x)$ forms an acute angle with the normal $\n(x)$ to $\Sigma$ at $x,$ not exceeding some constant which is
strictly smaller than $\pi/2.$  The construction of such field can be found, e.g., in \cite{HMT}, Proposition 2.3 (in fact a similar construction was used as long ago as in 1983 by A.P. Calder\'{o}n, see \cite{Calderon}).

For further reasoning, we need one more geometric construction. We suspect that it might have been well known since long ago, but were not able to locate a reference in the literature. Therefore, we describe it here.

\begin{lem}\label{lem.surface} Let $\Om\subset\R^{d+1}$ be a bounded domain with Lipschitz boundary $\Si.$ Let $\pmb{\G}$ be the transversal smooth vector field, as above.  In a neighborhood of $\Si$ denote by $\ell_{x}$ the integral curve of this vector field passing through $x$. Then, for a certain interior collar neighborhood $\Us$ of $\Si$ there exists a smooth surface $\Sib\subset \Us$, which has exactly one intersection point with every integral curve $\ell$ of the field $\pmb{\G}$ in $\Us,$ and this intersection is transversal.
\end{lem}
\begin{proof}
 In some two-sided collar neighborhood $\widetilde{\Us}$  of $\Si,$ the integral curves $\ell_{x'},$ $x'\in\Si,$ of the vector field $\pmb{\G}$ form a one-dimensional foliation. We define the function $\fF(x),$ $x\in \widetilde{\Us},$ as the distance of the point $x$ to the intersection point $x'\in\Si$ along the integral curve $\ell_{x'}$ passing through $x$ (with proper sign, positive inside $\Si.$)  Then $x=(x',t)$, $t=\fF(x)$ form a new co-ordinate system in $\widetilde{\Us}.$ In these co-ordinates the function $\fF$ is  continuous in $\widetilde{\Us}$, differentiable almost everywhere and has derivatives of any order in $t$ variable, in other words, along $\ell_{x'}$. Moreover, $\fF$ is strictly decaying in $t$ variable, i.e., along the curves $\ell_{x'}$.

 We take a smooth non-negative cut-off function $\chi_\e(x),\, \diam\supp\chi_\e<\e,$ $\int \chi_\e(x) dx=1,$ with sufficiently small $\e$, and consider the convolution, the smooth function $\fF_\e=\fF*\chi_\e.$ Sufficiently close to $\Si$, this function is smooth,  growing in $t$ variable, and  it has the same sign as $\fF$ outside the $\e$-neighborhood of $\Si$. We consider the level surface $\Sib$ for $\fF_\e,$ $\fF_\e(x)=\de$ for some $\de>0.$ By construction, being the level surface of a smooth function without stationary points, this surface is smooth, has only one intersection point with each of trajectories $\ell_{x'}$ and is uniformly nontangent with these trajectories. Due to the monotonicity of $\fF,$ the level $\de$ can be chosen in such way that  the surface $\Sib$ lies strictly inside $\Om;$ we may also suppose that $\dist(\Si, \Sib)>\rb_0>0.$
\end{proof}
\subsection{Change of variables}
We introduce now new co-ordinates in $\Om$ near $\Si$. Namely, for
the smooth surface $\Sib$, just constructed, for a point $x\in \widetilde{\Us},$ we take now  $x=(x', x_{d+1}),$ where $x'\in \Sib$  is the intersection point of the curve $\ell_x$ with $\Sib$ and $x_{d+1}$ is the distance from the point $x\in\Us$ to the surface $\Sib$ along the integral curve $\ell_x;$  we set $x_{d+1}$ to be positive outside $\Sib$.  According to the properties of the surface $\Sib,$ this is a smooth co-ordinates change from the initial Euclidean co-ordinates.

 In these co-ordinates, the initial surface $\Si$ is described \emph{globally} by the equation $x_{d+1}=\pmb{\psi}(x')$, $x'\in\Sib,$ with some Lipschitz function $\psb(x')>\rb_0>0$ for all $x'\in\Sib$ (recall that the surface $\Sib$ lies strictly inside $\overline{\Om}$).
 We denote by $\Cs_0$ the truncated cylinder with base $\Sib$ and top $\Si$:
  \begin{equation*}
  \Cs_0=\{(x',x_{d+1}): x'\in\Sib, x_{d+1}\in[0,\psb(x')]\}\subset{\Om},
  \end{equation*}
which is diffeomorphic to a collar neighborhood of $\Si$ in $\Om,$ and denote by $\Om_0$ the complement of this cylinder in $\Om$,
  \begin{equation*}
  \Om_0=\Om\setminus \Cs_0
  \end{equation*}


  In this setting the domain $\Om$ can be considered as a smooth manifold with Lipschitz boundary $\Si$,
   \begin{equation*}
    \Om =\Om_0\cup\Cs_0
   \end{equation*}
   having a truncated cylindrical exit $\Cs_0=\Cs_0[\psb]$  with base section $\Sib$ and  boundary $\Si$.


  Our next aim now is to perform a change of variables $\Fib$  of $\Om$ to the standard domain $\Omb=\Om_0\cup\Cs$ where $\Cs$ is the\emph{ straight cylinder}, $\Cs=[0,1]\times \Sib.$ This change of variables  is the identity mapping, $\Fib(x)=x$ on $\Om_0$ and an adjoining part of the cylinder $\Cs_0[\psb]$, $x_{d+1}<t_0$, while  near the Lipschitz boundary, it is glued together, by means of proper cut-offs, with the \emph{stretching} along the $x_{d+1}$ variable, i.e., along the integral curves of $\pmb{\G}$ namely, with the mapping
  \begin{equation}\label{change}
    \Fib(x)\equiv \Fib(x',x_{d+1})=(y',y_{d+1})\equiv\left(x', \frac{x_{d+1}}{\psb(x')}\right), \, x'\in\Sib,
  \end{equation}
 so,  at the top, $\Si$, we have $y_{d+1}=1$.
 Thus, in new co-ordinates, the Lipschitz boundary $\Si$ of the domain $\Om$ is transformed to the  smooth boundary of the domain $\Omb$: $y'\in\Sib, y_{d+1}=1.$

We denote by $x=\Psib(y)$  the mapping inverse to $\Fib$ (it, obviously, exists and it is the inverse stretching along the $y_{d+1}$ variable near $\Sib\times\{1\}$). The mapping $\Fib$ is is a bilipschitz mapping of $\Om$ to $\Omb$. Its restriction to $\Si$  is a bilipschitz mapping to $\Sib$. Such changes of variables leave invariant the Sobolev spaces $H^1$ and for the weak derivatives the usual chain rule holds, see, e.g., \cite{Ziemer}, Theorem 2.2.2.

The observation, crucial for our further reasoning is the following.

\begin{lem}\label{Lem.Change}
  \begin{enumerate}\item The components of the mappings $\Fib,$ $\Psib$ are continuous, their first order derivatives in $x'$, resp., $y'$   variables  and all their derivatives in ${x_{d+1}},$ resp., $y_{d+1}$ variables are bounded functions, smooth in $x_{d+1},$ resp., $y_{d+1}$ variable; they are continuous as functions of these variables with values in  $L_\infty(\Sib)$, resp. in $L_\infty(\Si).$
     \item The Jacobian matrices of the mappings $\Psib$ and $\Fib$ are bounded; their entries are continuous in $y_{d+1},$ resp., $x_{d+1},$ variables as functions with values in $L_\infty(\Sib)$, resp. in $L_\infty(\Si)$.
  \end{enumerate}
\end{lem}
All these properties follow automatically from the definition of the mapping $\Psib$ and the chain rule for derivatives, taking into account that the first order derivatives of a Lipschitz function are bounded.

 \section{Transformations and approximation of the P-S operator.}\label{Sect.NtD}
 In the course of the proof of our main result we will often transform operators under study while preserving their spectra, so that the initial and the approximating operator can be comfortably compared. Our considerations take place in the domain $\Om$ with cylindrical exit, as presented in the previous section. Having the quadratic forms $\rob_0$, $\ab_0$, we find out what happens with them under the transformation $\Psib.$

 \subsection{Transformation of quadratic forms}\label{Sect.trans}
 We study here how the quadratic forms defining the spectrum   of the Poincare-Steklov problem transform under the change of variables, as described in Sect. \ref{Sect.geom}; we set here  $v(y)=u(x),$ $y=\Fib(x)$, $x=\Psib(y)$, $x\in \Om,$  $y\in\Omb .$

 First, we consider the quadratic form
 \begin{equation*}
  \rob_0[u]=\rob[v]\equiv\int_{\Sib}\ro_0(\Psib(y))|v(y',1)|^2|\det J(y',1)|d\m_{\Sib}(y'),
 \end{equation*}
where $J(y',1)$ is the Jacobian matrix of the mapping $\Psib$ and $\m_{\Sib}$ is now the natural measure  on $ \Sib,$ a smooth surface in $\R^{d+1}$.
This quadratic form can be written as
\begin{equation*}
\rob[v]=   \int_{\Sib}{\ro}(y')|v(y',1)|^2d\mu_{\Sib}(y'),
\end{equation*}
with density ${\ro}(y')=\ro_0(\Psib(y',1))|\det J(y',1)|$.

 The quadratic form $\ab_0[u],$ see \eqref{form ab}, is transformed to
\begin{equation}\label{transform1}
  \ab_0[u]\equiv {\ab}[v]=\int_{\Omb}\sum_{j,k=1}^{d+1}\check{a}_{j,k}(y)\partial_j v(y)\overline{\partial_k{v(y)}}dy+
  \int_{\Omb} \check{\vb}(y)|v(y)|^2 dy,
\end{equation}
$v\in H^1(\Omb).$
The transformed coefficients $\check{a}_{j,k}(y), \check{\vb}(y)$, in the quadratic forms \eqref{transform1},  can be calculated, using \eqref{change}, but we do not need here their particular explicit expression at this point. What is important is the kind of dependence of these coefficients
on the variables $y_k$.

Recall that the coefficients $a_{j,k}(x)$ are bounded and are continuous at $\Si$ in the sense of \eqref{Def.Cont.bdry}. This property persists after the continuous change of variables $x=\Psib(y)$, so the functions $a_{j,k}(\Psib(y))$ are continuous at the new, smooth, boundary $\Sib$ in the sense of \eqref{Def.Cont.bdry}. Next, to obtain $\check{a}_{j,k}(y)$, we need to multiply $a_{j,k}(\Psib(y))$ by some algebraic combinations of derivatives of the mapping $\Psib,$ namely, the derivatives of $\Psib$, appearing in the chain rule when passing from $x$-derivatives to $y$-derivatives and also by the Jacobian of the mapping $\Psib$ arising in the passage from $dx$ in the integral to $dy$.  The  multiplication by these derivatives which are bounded but, generally, not continuous, destroys the continuity of the boundary values of $a_{j,k}(\Psib(y))$ for $y\in \Sib$. However, these derivatives are still continuous in $y_{d+1}$ variable as functions with values in $L_\infty(\Sib)$. We introduce the following definition.

\begin{defin}\label{Prop.noncont} Let $q(y)=q(y', y_{d+1})$ be a bounded measurable function in the variables $y'\in\Sib,$  $y_{d+1}\in(1-\de',1].$ Let ${\qb}(y')\in L_\infty(\Sib)$  be fixed. Denote by ${\qb^\diamond}(y)\equiv {\qb^\diamond}(y',y_{d+1})$ the continuation of ${\qb}(y')$ as a function not depending on $y_{d+1},$  ${\qb^\diamond}(y',y_{d+1}):={\qb}(y'),$ $y'\in\Sib.$ If for any $y'\in\Sib,$
\begin{equation}\label{def.infty.cont}
  \lim_{\de\to0}\|q-\qb^\diamond\|_{L_{\infty}(B((y',1),\de)}=0,
\end{equation}
we say that $q(y)$ has $L_\infty(\Sib)$ limit ${\qb}(y')$ at $\Sib$ in the variable $y_{d+1}\to 1$. If this is the case, we redefine the initial function $q(y)$ on $\Sib,$ by setting $q(y',1)=\qb(y')$, thus changing $q$ on a set of $d+1$-dimensional measure $0$. After this, we say that $q$ is $L_\infty(\Sib)$-continuous at $\Sib.$
\end{defin}
One can easily see that the product of  a function $L_\infty(\Sib)$-continuous at $\Sib$  and  a function continuous at $\Sib$ is again a function $L_\infty(\Sib)$-continuous at $\Sib.$
Therefore,  by Lemma \ref{Lem.Change}, the coefficients $\check{a}_{j,k}$ of the quadratic form $\ab[v]$ in \eqref{transform1} are $L_\infty(\Sib)$-continuous at $\Sib.$ This  'very week' continuity property proves to be sufficient for our approximation construction.

We recall that the ratio of quadratic forms $\rob_0/\ab_0$ defines an operator whose spectrum coincides with the one of the P-S operator. The domain of the quadratic forms is $H^1(\Omb)$. Therefore, our   transformations can be summarized in the following way.

\begin{proposition}\label{Prop.prop.coeff}
  The spectrum of the variational  Poincar\'e-Steklov problem  \eqref{LSteklov} coincides with the spectrum of the ratio
\begin{equation}\label{ratio5}
\frac{{\rob}[v]}{{\ab}[v]}, \, v\in H^1(\Omb),
\end{equation}
in a domain $\Omb$ with smooth boundary $\Sib$, and where the quadratic form $\ab[v]$ is elliptic and has bounded  coefficients, which are $L_\infty(\Sib)$ - continuous in $y_{d+1}$ variable at the boundary $\Sib=\partial\Omb.$
\end{proposition}

From now on, we may forget about the initial Lipschitz problem and study, from scratch, the eigenvalue distribution for the ratio \eqref{ratio5}. To simplify notations, we drop further on the 'check' over the symbols and, again, denote the main domain by $\Om$ and its, now smooth, boundary by $\Si$ ($=\Si\times\{1\}$ in our $y$ co-ordinates). Recall that we have paid for the smoothness of the boundary by weakening  the continuity property    of the coefficient matrix at the boundary, from continuity to $L_\infty(\Si)$ continuity. In the notations of the quadratic forms $\rob_0,$ $\ab_0$ we drop the subscript '0'.

\subsection{An operator representation of the eigenvalue problem.}

We recall that we  consider the case of a non-negative  function $\ro.$

The quadratic form $\ab[v]$ with domain $H^1(\Om)$ defines the self-adjoint positive 'Neumann' operator $\Tb=\Tb_\aF$ in $L_2(\Om).$ This elliptic operator acts, formally, as
 \begin{equation*}
  \Tb v\equiv\Tb_\aF v=-\sum_{j,k}\partial_j a_{j,k}\partial_k v +\vb v,
 \end{equation*}
in the sense of distributions, but its domain is rather hard to describe due to the lack of  continuity of the coefficients. Fortunately, neither the domain nor the exact action are  needed in our reasoning.
Instead, since $\Tb$ is a positive self-adjoint elliptic operator with bounded coefficients, the domain of the square root $\Tb^{\frac12}$ can be described explicitly: it coincides with $H^1(\Om)$ and $\Tb^{\frac12}$ is an isomorphism from $H^1(\Om)$ to $L_2(\Om)$.

It was already stated that we are free in the choice of the function $\vb$. Now we fix this choice, supposing that $\vb\in C_0^{\infty}(\Om)$ (after the above co-ordinates change).

Since $\ab[v]=(\Tb^{\frac12}v,\Tb^{\frac12}v),$ we can set $f=\Tb^{\frac12}v\in L_2(\Om)$ and therefore the ratio \eqref{ratio5}
transforms to
\begin{equation}\label{ratio1}
\frac{\rob[\Tb^{-\frac12}f]}{\|f\|_{L_2(\Om)}^2}=\frac{\int_{\Si}|(\g \Tb^{-\frac12}f)(y')|^2 \ro(y')d\m_{\Si}(y')}{\|f\|_{L_2(\Om)}^2}, \,f\in L_2(\Om).
\end{equation}
Here and further on, we denote by $\g:f\mapsto \g f$ the trace operator acting as the restriction to $\Si$ of a function $f$ defined on the domain  $\Om$ with cylindric exit $\Cs$ having  the boundary lid $\Si$ identified with $\Si\times\{1\}$. We recall that the boundary $\Si$  is smooth. Therefore, this restriction is known to act as a bounded operator from $H^s(\Om)$ to $H^{s-1/2}(\Si)$ for all \emph{positive} $s $ \emph{except the half-integer ones} ($s\in \N+\frac12).$ In particular, this implies that the operator $\g\Tb^{-\frac12}$ is bounded as acting from $L_2(\Om)$ to $L_2(\Si)$ (and even to $H^{\frac12}(\Si)$). Following Sect.2, we restrict ourselves to  $\ro\ge 0$ and $\ro\in L_\infty,$ therefore, the operator $\ro^{\frac12}\g\Tb^{-\frac12}: L_2(\Om)\to L_2(\Si)$ is bounded as well.

As a result, the numerator in \eqref{ratio1} can be written as
\begin{gather*}
{\rob[\Tb^{-\frac12}f]}=(\ro(y')^{\frac12}\g \Tb^{-\frac12}f,\ro(y')^{\frac12}\g \Tb^{-\frac12}f)_{L_2(\Si)}=\\\nonumber\left(\left[\ro^{\frac12}\g\Tb^{-\frac12}\right]^* [\ro^{\frac12}\g\Tb^{-\frac12}]f,f\right)_{L_2(\Om)}.
  \end{gather*}
In this way, we have reduced our spectral problem to the study of the spectrum of the self-adjoint operator
\begin{equation*}
\Gb_\aF=\left[\ro^{\frac12}\g\Tb^{-\frac12}\right]^*\left[ \ro^{\frac12}\g\Tb^{-\frac12}\right],
\end{equation*}
acting in $L_2(\Om)$.
Since the nonzero discrete spectrum of the product does not change under the cyclic permutation of the factors,
the nonzero eigenvalues  of $\Gb_\aF$ are the same as the ones of
the self-adjoint operator
\begin{equation}\label{Permuted}
  \Hb_\aF\equiv \left[\ro^{\frac12}\g\Tb^{-\frac12}\right] \left[\ro^{\frac12}\g\Tb^{-\frac12}\right]^*,
\end{equation}
acting in $L_2(\Si).$

It is more convenient to have a somewhat different representation of the operator $\Hb_\aF.$
\begin{lem}\label{lem.transpos}
  The following equality holds
  \begin{equation}\label{trans}
  \Hb_\aF=  \left[\ro^{\frac12}\g\Tb^{-\frac12} \right] \left[\ro^{\frac12}\g\Tb^{-\frac12}\right]^*=\ro^{\frac12}\g(\g\Tb^{-1})^* \ro^{\frac12}=\ro^{\frac12}\g(\ro^{\frac12}\g\Tb^{-1})^* .
  \end{equation}
\end{lem}
\begin{proof}We start with the equality
\begin{equation}\label{Adj.3}
[\ro^{\frac12}\g \Tb^{-\frac12}] \Tb^{-\frac12}= \ro^{\frac12}\g \Tb^{-1}: L_2(\Om)\to L_2(\Si),
\end{equation}
where we have the product of two bounded operators on the left. Next we take adjoint of both parts in \eqref{Adj.3},
\begin{equation}\label{Adj.4}
  \Tb^{-\frac12}[\ro^{\frac12}\g \Tb^{-\frac12}]^*=(\ro^{\frac12}\g \Tb^{-1})^*:L_2(\Si)\to L_2(\Om).
\end{equation}
This is a bounded operator acting from $L_2(\Si)$ to $H^1(\Om)$, again a product of two bounded operators. Therefore, the restriction operator $\g$ is well defined on the range of the operator in \eqref{Adj.4} and we can apply $\g$ on the left, and then multiply by the bounded function $\ro^{\frac12},$
which gives us \eqref{trans}.
\end{proof}

\subsection{Smooth approximation}\label{approximation}

Our study of P-S eigenvalues will be based upon a special smooth approximation of the coefficient matrix $\aF.$
\begin{lem}\label{lem.appr} Let $\aF\in L_\infty(\Om)$ be a Hermitian matrix-function, continuous at $\Si$ as a function of $y_{d+1}$ variable with values in $L_\infty(\Si)$ in the sense of Definition \ref{def.infty.cont}.   Then
  for any $\e>0$ and $p<\infty, $ there exist a  $C^\infty(\overline{\Om})$-smooth matrix-function $\tilde{\aF}(y),$ $y\in\bar{\Om},$ such that  the following conditions are satisfied
  \begin{enumerate}
  \item $\|\aF(.,1)-\ah(.,1)\|_{L_p(\Si)}<\e$;
  \item the function $\aF(y)-\ah(y)$  is continuous at $\Si$ as a function of $y_{d+1}$ variable with values in $L_\infty(\Si);$
  \item $\|\ah\|_{L_{\infty}(\Om)}\le  c\|\aF\|_{L_{\infty}(\Om)}$;
  \item if the matrix  $\aF$ is uniformly positive in $\Om,$ then, additionally, $\|\ah^{-1}\|_{L_{\infty}(\Om)}\le  c\|\aF^{-1}\|_{L_{\infty}(\Om)}.$
  \end{enumerate}
  with the constant $c$ not depending on $\aF.$
    \end{lem}
  The condition (3) means that    the $L_\infty$ norm for the approximating matrix can be arranged to be not depending on $\e$, being controlled, by the $L_\infty$ norm of the initial matrix $\aF;$  by condition (4), in the elliptic case, the same can be arranged for  the ellipticity constant of the approximating matrix.  Note that the approximation quality in Lemma \ref{lem.appr} is required only at the boundary $\Si$.

  \begin{proof} We discuss the ellipticity case, $\aF(y)\ge C>0,$ first. Consider the constant matrix $\aF_1(y)=C_1\mathbf{I}_{(d+1)\times(d+1)}$,
  where $\|\aF^{-1}\|_{L_\infty(\Om)}^{-1}\le C_1\le \|\aF\|_{L_\infty(\Om)}.$ Next, to construct $\ah(y)$ near the boundary, we take $\de>0$ such that
    for any ball $B(y',\de), \, y'\in\Si$ the  inequality
  $\|\aF (y)-\aF^\diamond(y)\|_{L_\infty(B(y',\de))}<\e$ holds (where ${\aF}^\diamond(y)$ is the continuation of $\aF(y',1)$ from $\Si$ to the $\delta$-neighborhood of $\Si,$
${\aF}^\diamond(y',y_{d+1})=\aF(y',1).$)
Next we  approximate the matrix $\aF(.,1)$ by a smooth matrix $\ah_0\in L_p(\Si)\cap L_{\infty}(\Si)$
 such that $\|\ah_0(.)-\aF(.,1)\|_{L_p(\Si)}<\e,$ thus the condition 1)  in Lemma \ref{lem.appr} is satisfied. Also we require that
   \begin{equation}\label{twosided}
   \|\ah_0\|_{L_\infty(\Si)}\le C\|\aF\|_{L_\infty(\Om)}, \, \|\ah_0^{-1}\|_{L_\infty(\Si)}\le C'\|\aF^{-1}\|_{L_\infty(\Om)};
 \end{equation}
all these conditions can be satisfied, for example, by setting $\ah_0(.) = \Vb(\tau)\aF(.,1)$, where $\Vb(\tau)$ is the heat semigroup on $\Si,$  since $\Vb(\tau)\to \pmb{1}$ strongly in $L_p(\Si),$ $p<\infty$  as $\tau\to +0$ and $\Vb(\tau)$ is bounded in $L_\infty,$ $\tau\ge 0.$
  We consider now the matrix ${\ah_0}^\diamond(y),$ the continuation of $\ah_0(y'),$ namely ${\ah_0}^\diamond(y', y_{d+1})=\ah_0(y').$ This smooth matrix, not depending on $y_{d+1}$ is, of course, continuous in the $y_{d+1}$ variable  as a function with values in $L_\infty(\Si)$. Therefore, the difference, $\aF(y)-{\ah_0}^\diamond(y)$ is also continuous at $\Si$ as a function  of $y_{d+1}$  with values in $L_\infty(\Si)$. This grants the condition 2) in Lemma.  Conditions 3) and 4) are satisfied due to \eqref{twosided}. Finally, we glue together the matrix ${\ah_0}^\diamond$ defined in a neighborhood of $\Si$ and the matrix $\aF_1(y)$  defined outside such neighborhood by means of proper cut-off functions, and this matrix satisfies all conditions of Lemma.

  In the case when the matrix $\aF$ is not supposed to be uniformly positive definite, the simplified reasoning succeeds: namely we just skip all estimates concerning $\aF^{-1}$, therefore, we may simply set $\aF_1(y)=0$.
  \end{proof}

We will use Lemma \ref{lem.appr} twice; first time when approximating the elliptic coefficients  matrix $\aF(y),$ where the control of the ellipticity constant is important, and the second time for the approximation of the matrix $\bFb=\aF-\ah$ when proving the required operator  convergence of $\Hb_{\ah}$ to $\Hb_{\aF}$.

\subsection{The approximating operator}
Consider, for a given $\e>0$, the approximating coefficient matrix $\tilde{\aF}$ constructed in Section 4.3 for the elliptic coefficient matrix $\aF$ and the corresponding elliptic differential operator $\Lch=\Lc_{\tilde{\aF}},$ \emph{with the same zero order term} $\vb\in C^{\infty}(\overline{\Om}).$
With the coefficient  matrix $\tilde{\aF}$ we associate the Neumann operator $\tilde{\Tb}=\Tb_{\tilde{\aF}},$ by means of the quadratic form $\pmb{\ah}[v],$ $v\in H^1(\Om)$.

Since the boundary of $\Om$ is  smooth,   as well as the   coefficients matrix  $\ah$, some more can be said about the properties of the operator \eqref{Permuted} for this matrix.
The corresponding compact operator, which we denote by ${\Hb}_{\ah},$ has, by Lemma \ref{lem.transpos}, the form
\begin{equation}\label{tildeH}
   {\Hb}_{\ah}=\ro^{\frac12}\g(\g\Tb_{\ah}^{-1})^* \ro^{\frac12}.
\end{equation}
First of all, for the eigenvalues of this operator, which coincide with eigenvalues of the corresponding $\ND$ operator, the asymptotic formula holds, of the form  \eqref{formula}, \eqref{beta}, \eqref{alpha}, with the matrix $\ah$ replacing $\aF,$ due to the existing results, see Theorem \ref{Th.Agr}.
Next, since the coefficients of  $\Tb_{\ah}$ are smooth, the domain and action of the operator $\Tb_{\ah}$ can be described explicitly. Namely, by the standard elliptic regularity results,  the domain $\Dc(\tilde{\Tb})$ of $\tilde{\Tb}$ is the subspace in the Sobolev space $H^2(\Om)$, consisting of functions $v$ satisfying the classical Neumann boundary conditions, $\sum_k \tilde{a}_{j,k}\n_k\partial_j v|_{\partial \Om} =0,$ where $\n_k$ are the components of the unit  normal vector to $\Si=\partial\Om.$

\begin{lem}\label{Lem2} For a smooth matrix $\ah,$ for any $f\in L_2(\Si)$, the function $v=(\ro^{\frac12}\g\tilde{\Tb}{}^{-1})^*\ro^{\frac12} f$ satisfies in $\Om$
the elliptic equation
\begin{equation}\label{smooth}
  \tilde{\Lc}v\equiv\tilde{\Lc}  (\ro^{\frac12}\g\tilde{\Tb}{}^{-1})^*\ro^{\frac12}f =0 \,\, \mathrm{in}\,\, \Om.
\end{equation}
\end{lem}
\begin{proof}Let $g\in C^\infty_0(\Om)$ be a smooth function with compact support in $\Om.$ Consider  the sesquilinear form $\Ib[f,g]=(\tilde{\Lc}(\ro^{\frac12}\g\tilde{\Tb}{}^{-1})^*\ro^{\frac12}f,g)_{L_2(\Om)}.$ The function $g$ belongs to the domain of the self-adjoint operator $\tilde{\Tb},$ therefore,
\begin{equation}\label{homog}
  \Ib[f,g]=((\ro^{\frac12}\g\tilde{\Tb}{}^{-1})^*\ro^{\frac12}f, \tilde{\Tb}{}g)_{L_2(\Om)}= (\ro^{\frac12}f,\ro^{\frac12}\g\tilde{\Tb}{}^{-1}
  \tilde{\Tb g} )_{L_2(\Si)}=(f,\ro\g g)_{L_2(\Si)}.
  \end{equation}
Since $g=0$ on $\Si,$ the last expression in \eqref{homog} is zero. Therefore, the function $\tilde{\Lc}(\ro^{\frac12}\g\tilde{\Tb}{}^{-1})^*\ro^{\frac12}f$ is orthogonal to $C^\infty_0(\Om),$ and this implies that $\tilde{\Lc}(\ro^{\frac12}\g\tilde{\Tb}{}^{-1})^*\ro^{\frac12}f$ is zero in $\Om.$
\end{proof}

\subsection{$\Hb_{\aF}$ as a pseudodifferential operator}
Here we establish an important property of the approximating operator $\Hb_{\tilde{\aF}},$ needed further for evaluating the rate of approximation of spectra.

\begin{proposition}\label{Prop.Pseudodiff}
For a smooth approximating matrix $\tilde{\aF}$ and a smooth weight $\ro$, the operator $\Hb_{\tilde{\aF}}$ is an  order $-1$ pseudodifferential operator on the boundary $\Si.$
\end{proposition}
\begin{proof}
 When finding, in the smooth case, an expression for the operator $\Hb_{\tilde{\aF}},$ we study its component $(\ro^{\frac12}\g \tilde{\Tb}^{-1})^*$ first.
 Recall that $\tilde{\Tb}{}^{-1}$ is the inverse of the realization of the operator $\Lch$ with Neumann boundary conditions, $\g_1 u=0$, where $\g_1$ is the conormal derivative. Continue the coefficients $\ah$ outside $\Om$ in a smooth way.  We denote by $\Rb$ the fundamental solution for $\Lch$; so that $\Lch \Rb-1,$ $\Rb\Lch-1$ are  infinitely smoothing operators; $\Rb$ is the order $-2$ pseudodifferential operator with symbol $\rb(y,\y)$ having the principal term $\rb_{-2}(y,\y)=\pmb{\ah}_y(\y)^{-1},$ where $\pmb{\ah}_y(\y)=\sum_{j,k}\tilde{a}_{j,k}(y)\y_j\y_k$ is the principal symbol of the operator $\Lch$. Let $d>1$ first. Then the operator $\Rb$ acts on functions on $\Om$ as an integral operator with kernel $R(y,z-y) $ having leading singularity of order $1-d$ at the  diagonal $y=z$. Namely, the principal term,   $R_{1-d}(y,y-z),$  the  leading singularity in $R(y,z-y) $ is the Fourier transform of the principal symbol $\rb_{-2}(y,\y)$,
  \begin{equation*}
  R_{1-d}(y, y-z)=\Fc_{\y\to y-z}\rb_{-2}(y,\y).
  \end{equation*}
  For the dimension $d+1=2$, the kernel $R(y,z-y)$ has logarithmic singularity at the diagonal.

 Let $f$ be a function in $L_2(\Om)$. In order to find a representation for $u=\tilde{\Tb}{}^{-1}f,$ we set $u_0=\int_{\Om}R(y, y-z)f(z)dz$. This function satisfies the equation $\Lch u_0=f$ in $\Om$ but not the Neumann boundary condition $\g_1 u=0$, where $\g_1$ is the conormal derivative at the boundary, corresponding to the elliptic operator $\Lch$. We construct the correction $u_1,$ such that $\Lch u_1=0,$ $\g_1u_1=-\g_1u_0$. To do this,  we consider the \emph{Green function} $\Gc_1(y,z)$ for the Neumann problem, in other words, the integral kernel of the solution operator $\pmb{G}_1$ of the Neumann problem for $\Lch$ . This means that for a smooth function $h$ on $\Si$,
 \begin{equation}\label{Compens}
 \Lch\pmb{G}_1 h:=  \Lch_y\int_{\Si} \Gc_1(y,z)h(z)dz =0, \, y\in\Om,
 \end{equation}
 and
\begin{equation}\label{Green.Neu}
 \lim_{y\to y_0\in\Si}\g_{1,y}\int_{\Si} \Gc_1(y,z)h(z)dz=h(y_0).
\end{equation}
   The operator $\pmb{G}_1:C^\infty(\Sigma)\to C^\infty(\Omega)$ is a Poisson operator in the Boutet-de-Monvel algebra.
 By \eqref{Green.Neu},  it satisfies $\g_1\pmb{G}_1:C^\infty(\Si)\to C^{\infty}(\Si)=\pmb{1},$ i.e., it gives the identity operator.
 Therefore the solution of the boundary problem $\Lch u=f$ in $\Omega,$ $\g_1 u =0$ on $\Si$ can be expressed as
 \begin{gather*}
   u(y)=u_0(y)+u_1(y)=\\\nonumber \int_{\Om}R(y,y-z) f(z)dz- \int_{\Si} \Gc_1(y,z)\g_1\left(\int_{\Om}R(y,y-z) f(z)dz\right )\\\nonumber
   =(\Rb f)(y)-(\pmb{G}_1\g_1\Rb f)(y).
 \end{gather*}
 As a result, the operator $\tilde{\Tb}{}^{-1}$ can be represented as
  \begin{equation*}
   \tilde{\Tb}{}^{-1}=\Rb-\pmb{G}_1\g_1\Rb=(1-\pmb{G}_1\g_1)\Rb;
  \end{equation*}
 such representation,  in a somewhat different setting, can be found, e.g.,  in   \cite{Grubb.book}, Theorem 9.20. There, after the standard straightening of the boundary,  the operator $\Rb$ is treated as a truncated pseudodifferential operator,
 \begin{equation}\label{psdo0}
   \Rb=\chi_\Om \RF \eb_{\Om},
 \end{equation}
where $\eb_{\Om}$ is the operator of extension by zero of a function in $\Om$ to the whole $\R^{d+1}$, $\RF$ is a pseudodifferential operator in $\R^{d+1}$ with symbol $\rb(y, \y)$ and $\chi_\Om$ is the restriction of functions in $\R^{d+1}$ to $\Om.$

We also need an expression for the operator $\tilde{\Tb}_0^{-1}$, the resolvent of  the \emph{Dirichlet} problem for the equation $\Lch u=f$ in $\Om.$ In a similar way,
\begin{equation}\label{Yb6}
  \tilde{\Tb_0}{}^{-1}=\Rb-\pmb{G}_0\g\Rb=(1-\pmb{G}_0 \g)\Rb,
  \end{equation}
where $\pmb{G}_0$ is the Poisson operator (the Green function) solving  the nonhomogeneous Dirichlet boundary problem for $\Lch$ in $\Om$:
$\Lch u=0,\, \g u=h;$ this means that  $\g \pmb{G}_0=\pmb{1}.$

We pass now to $\g \tilde{\Tb}{}^{-1}.$ Since $\g {\tilde{\Tb}_0{}^{-1}}=0,$
we can write

\begin{gather}\label{Yb7}
  \g\tilde{\Tb}{}^{-1}=\g(\tilde{\Tb}{}^{-1}-\tilde{\Tb}_0{}^{-1})=\g(\pmb{G}_0\g\Rb-\pmb{G}_1\g_1\Rb)\\\nonumber =\g\pmb{G}_0\g\Rb-\g\pmb{G}_1\g_1\Rb=
  (\g\pmb{G}_0\g-\g\pmb{G}_1\g_1)\Rb.
\end{gather}
We recall now that $\g\pmb{G}_0=\pmb{1}.$ Further on, the function $u_1=\pmb{G}_1\g_1\Rb f$ satisfies the equation $\Lch u_1=0$ in $\Om$, therefore $\g u_1$ and $\g_1 u_1$ are connected by the Neumann-to Dirichlet operator $\ND$, $\g u_1=\ND\g_1 u_1$. We set all this into \eqref{Yb7} and obtain
\begin{equation}\label{Yb8}
  \g\tilde{\Tb}^{-1}=(\g-\ND\g_1)\Rb.
\end{equation}
Consequently,
\begin{equation*}
  (\g \tilde{\Tb}^{-1})^*= (\g\Rb)^* -(\ND\g_1 \Rb)^*. 
\end{equation*}
Finally, we have
\begin{equation}\label{Yb91}
  \g (\g \tilde{\Tb}{}^{-1})^*=\g(\g\Rb)^*-\g(\ND\g_1 \Rb)^*.
\end{equation}
Both terms on the right in \eqref{Yb91} belong to the Boutet-de-Monvel algebra, see e.g.,  \cite{Grubb.book}, where the composition rules are described in detail. In this setting,  $\Rb$ is a truncated pseudodifferential operator of order $-2$, $\g$ is a trace operator of zero order, further, $(\g\Rb)^*$ is a Poisson operator, and as a  result,  $\g(\g\Rb)^*$ is a pseudodifferential operator  of order $-1.$ Similarly,  $\g_1$ is a trace operator of order $1$, the Neumann-to Dirichlet operator $\ND$ is an order $-1$ pseudodifferential operator  on the boundary, therefore, the second term in \eqref{Yb91} is a pseudodifferential operator of order $-1$ as well. Its principal symbol can be expressed in a standard way, algebraically, via the principal symbols of $\Rb$ and $\ND$, as we will see in Sect.7. Finally, the multiplication by the smooth function $\ro^{\frac12}$ produces, again, a pseudodifferential operator.
\end{proof}

  \section{Operator perturbations}\label{Sect.Pert.}
  Our approach to establishing the asymptotic formula for Steklov eigenvalues, similarly to \cite{RT1}, is based upon the operator approximation.
  This approximation idea was successfully used by M.S. Agranovich in  \cite{Agr06} on the base of the variational setting of the problem. In fact,  for  more smooth boundaries, the ones of class $C^{1,1}$, a similar  perturbational approach was even  used in  \cite{Sandgren}, the earliest paper on the Steklov eigenvalue asymptotics.

   Namely, if the coefficients matrix $\aF(y)$ is continuous, it can be approximated in $C(\Om)$, both from above and from below, by smooth elliptic  matrices $\ah_{\pm}(y)$, and the required closeness of spectra of operators $\Tb$ and $\tilde{\Tb}$ follows in \cite{Sandgren} and in \cite{Agr06} from rather simple monotonicity estimates. The weakening of conditions imposed on the boundary in \emph{this} paper, actually, the passage from the \emph{ continuous} derivatives of the function $\psi$ defining the boundary to the function having only \emph{bounded} derivatives, seemingly, a minor one, is, in fact, rather essential, the resulting coefficients $a_{j,k}$ are not continuous any more and they cannot be approximated in the $C$ metric by smooth ones. Instead, we use an approximation in a weaker, $L_p$ sense, as in Lemma \ref{lem.appr}, which turns out to be sufficient.

  \subsection{Approximation of the operator}

  In this section, we find the expression for the difference of operators describing the eigenvalues of the N-to-D operators. Let $\aF,\ah$ be the matrices of coefficients of the operators $\Lc, \tilde{\Lc}$, described in Sect.4, so that $\aF,\ah{}^{-1}$ belong to $L_\infty(\Om)$, $\ah,\ah^{-1}\in C^\infty(\overline{\Om})$  and $\ah|_{\Si}- \aF|_\Si$ is small in the sense of Lemma \ref{lem.appr}.

   Consider the operators $\Tb, \tilde{\Tb}$, the Neumann operators for $\Lc, \Lch$.
  \begin{lem}\label{LemDiffRes} Under  the above conditions, the following factorization is valid:
   \begin{equation}\label{factorization}
    \Tb^{-1}-\tilde{\Tb}{}^{-1}= \Xb^*\Yb,
   \end{equation}
where $\Xb,\Yb$ are bounded operators acting from $L_2(\Om)$ to $L_2(\Om)\otimes \C^{d+1}$.
  \begin{equation}\label{factors}
    \Xb=\nabla\Tb^{-1}, \Yb= \bFb(y)\nabla\tilde{\Tb}{}^{-1},
  \end{equation}
  with $\bFb=(\ah-\aF).$
 \end{lem}
 \begin{proof}
   Essentially, our Lemma \ref{LemDiffRes} is an  analogy of Lemma 8.1 in \cite{BS}. For the sake of completeness, we present the detailed  proof of \eqref{factorization} in our setting.

   The equality \eqref{factorization} is equivalent to
   \begin{equation}\label{factor.1}
     (\Tb^{-1} f,g)_{L_2(\Omega)}-(\tilde{\Tb}{}^{-1}f,g)_{L_2(\Omega)}= \int_{\Om}\al (\ah-\aF)\nabla\tilde{\Tb}{}^{-1}f,\nabla\Tb^{-1}g \ar dy, \, f,g\in L_2(\Om),
   \end{equation}
   where the angle brackets denote the scalar product in $\C^{d+1}.$

    In our conditions, the domains of the operators $\Tb^{\frac12}, \tilde{\Tb}{}^{\frac12}$ coincide, both are equal to the Sobolev space $H^1(\Om).$ It follows, in particular, that $\Dc(\tilde{\Tb})\subset H^2(\Om)\subset \Dc(\Tb^{\frac12}).$ Let $u,v$ be arbitrary functions in $\Dc(\Tb^{\frac12})=\Dc(\tilde{\Tb}{}^{\frac12})=H^{1}(\Om).$ Consider the equality
    \begin{gather}\label{factor.2}
      (\tilde{\Tb}{}^{\frac12}u,\tilde{\Tb}{}^{\frac12}v)_{L_2(\Om)}- (\Tb^{\frac12}u,\Tb^{\frac12}v)_{L_2(\Omega)}=\tilde{\ab}[u,v]-\ab[u,v]=\\\nonumber
      \int_\Om \al \ah\nabla u,\nabla v\ar dy-
      \int_\Om \al \aF\nabla u,\nabla v\ar dy\\\nonumber  =\int_{\Om}\al(\ah-\aF)\nabla u,
      \nabla v\ar dy.
    \end{gather}
   We set here $u=\tilde{\Tb}{}^{-1}f,$ $v=\Tb^{-1}g,$ where $f,g$ are arbitrary elements in $L_2(\Om)$.  These functions $u,v$ belong to $H^1(\Om)$, therefore \eqref{factor.2} is satisfied. Such substitution leads to \eqref{factor.1}:
   \begin{equation*}
     (\tilde{\Tb}^{\frac12}u,\tilde{\Tb}^{\frac12}v)_{L_2(\Omega)}=
     (\tilde{\Tb}^{\frac12}\tilde{\Tb}^{-1}f,\tilde{\Tb}^{\frac12}{\Tb}^{-1}g )_{L_2(\Omega)}=(f,{\Tb}^{-1}g )_{L_2(\Omega)}=({\Tb}^{-1}f,g)_{L_2(\Omega)},
   \end{equation*}
   and, similarly, for $(\Tb^{\frac12}u,\Tb^{\frac12}v),$
   \begin{gather*}
     ({\Tb}^{\frac12}u,\tilde{\Tb}^{\frac12}v)_{L_2(\Omega)}=(\Tb^{\frac12}\tilde{\Tb}^{-1}f,\Tb^{\frac12}{\Tb}^{-1}g)_{L_2(\Omega)}=\\\nonumber
     (\Tb^{\frac12}\tilde{\Tb}^{-1}f,\Tb^{-\frac12}g )_{L_2(\Omega)}=(\tilde{\Tb}^{-1}f,g )_{L_2(\Omega)}.
   \end{gather*}
 \end{proof}

 Using \eqref{factors}, we arrive at the representation of the difference of the operators $\Hb_a$ and $\Hb_{\ah}$, see \eqref{Permuted}, \eqref{tildeH}:
 \begin{equation}\label{factor.3}
   \Hb_\aF-\Hb_{\ah} \equiv\ro^{\frac12}\g( \g(\Tb^{-1}-\tilde{\Tb}{}^{-1}))^*\ro^{\frac12}=\ro^{\frac12} \g[\g\Xb^* \Yb]^*\ro^{\frac12}.
 \end{equation}

 We need a somewhat different representation for the operator in \eqref{factor.3}.

 \begin{lem}\label{lem.factor}The following equality is valid

  \begin{equation}\label{factor.4}
   \Hb_\aF-\Hb_{\ah}=\ro^{\frac12}(\g \Yb^*)(\g \Xb^*)^*\ro^{\frac12}
. \end{equation}
 \end{lem}
 \begin{proof}
   The operator $\g\Xb^*\Yb =(\g\Xb^*)(\Yb)$ is a product of two bounded operators, therefore, $(\g\Xb^*\Yb)^*=\Yb^*(\g\Xb^*)^*.$ We apply the operator $\g$ from the left, multiply by the bounded function $\ro^{\frac12},$  and obtain \eqref{factor.4}.
 \end{proof}

We will use  expression \eqref{factor.4} to evaluate the singular numbers of the difference $\Hb_\aF-\Hb_{\ah}.$ Since the function $\ro$ is bounded, the multiplication by $\ro^{\frac12}$ preserves spectral estimates, therefore, it suffices to drop this weight in further estimates.
 \subsection{Spectral estimates for the operator $\g\Xb^*$}

 We represent the operator $\Xb=\nabla \Tb^{-1}$ in the following way as a product of two bounded operators,
 \begin{equation*}
   \Xb=(\nabla \Tb^{-\frac12} )\Tb^{-\frac12}, \Xb^*=\Tb^{-\frac12}(\nabla \Tb^{-\frac12})^*.
 \end{equation*}
 This gives
 \begin{equation*}
   \g\Xb^*=(\g\Tb^{-\frac12})(\nabla \Tb^{-\frac12})^*.
 \end{equation*}
 The operator $(\nabla \Tb^{-\frac12})$ is bounded in $L_2(\Om)$ since $\Tb^{-\frac12}: L_2(\Om)\to H^1(\Om)$ and $\nabla$ is bounded as acting from $H^1(\Om)$ to $L_2(\Om).$ The norm of $\nabla \Tb^{-\frac12}$  is controlled by the ellipticity constant of $\Lc$ and  does not depend on the approximation $\ah$. Therefore, $s$-numbers of the operator $\g \Xb^*$ are majorated by $s$-numbers  of $\g\Tb^{-\frac12},$
\begin{equation*}
\nb^{\sup}(2d, \g \Xb^*)\le C \nb^{\sup}(2d, \g\Tb^{-\frac12} ).
\end{equation*}
 Next we have
 \begin{equation*}
   n(\la,\g\Tb^{-\frac12})=n(\la^2, (\g\Tb^{-\frac12})^*(\g\Tb^{-\frac12})).
 \end{equation*}
 Therefore,
 \begin{equation}\label{X3}
   \nb^{\sup}(2d, \g\Tb^{-\frac12})=\nb^{\sup}(d,(\g\Tb^{-\frac12})^*(\g\Tb^{-\frac12}))
 \end{equation}
 The expression on the right in \eqref{X3} characterizes the counting function for the eigenvalues of the  ratio
 \begin{equation}\label{X4}
   \frac{\|\g \Tb^{-\frac12} f\|^2_{L_2(\Si)}}{\|f\|^{2}_{L_2(\Om)}}, \, f\in L_2(\Om).
 \end{equation}
  We set $f=\Tb^{\frac12}u$, $u\in H^1(\Om) $ in \eqref{X4} and obtain the  ratio
  \begin{equation}\label{X5}
    \frac{\|\g u\|^2_{L_2(\Si)}}{\|\Tb^{\frac12}u\|_{L_2(\Om)}}=\frac{\int_{\Si}|u(y',1)|^2d\m_{\Si}(y')}{\ab[u]}\le \frac{C^{-1}\int_{\Si}|u(y',1)|^2d\m_{\Si}(y')}{\|u\|^2_{H^1(\Om)}}, \, u\in H^1(\Om),
  \end{equation}
 since $\ab[u]\ge C\|u\|^2_{H^1(\Om)}$, with constant determined by ellipticity constant of $\Lc,$ the ratio \eqref{X5} is majorated by  the spectral problem considered in Theorem \ref{Th.estimate}. From this theorem, taking into account \eqref{X3}, we obtain
 \begin{equation}\label{X6}
  \nb^{\sup}(2d, \g\Tb^{-\frac12})\le C,
 \end{equation}
 with constant $C$ depending only on the $L_\infty$ norm of the matrix $\aF^{-1},$ therefore, returning to the original problem, only on the ellipticity constant of the initial operator and the Lipshitz norm of the function $\psi.$

 \section{Spectral estimates for the operator $\g\Yb^*$. Rough estimates}\label{Estimates Yb}
   In the study of the second term in the factorization \eqref{factor.3}, \eqref{factor.4}, namely,  $\g\Yb[\bFb]^*$, where we set $\Yb[\bFb]=\bFb\nabla \tilde{\Tb}{}^{-1},$ we will need to obtain an asymptotic singular numbers  estimate for $\nb^{\sup}(2d,\g\Yb^*)=\limsup_{\la\to 0}{\la^{2d}n(\la, \g\Yb^*)}.$ This will be done in two steps. First, in this section, for an \emph{arbitrary} matrix function $\bFb$, $L_\infty(\Si)$-continuous at $\Si\times\{1\},$ we find an estimate of $\nb^{\sup}(2d,\g\Yb^*)$ in terms of certain integral norm of $\bFb(.,1).$ The constant in this estimate will, unfortunately,  depend, in an uncontrollable manner, on the approximating matrix $\ah$.  However, this estimate  enables us to restrict ourselves to smooth matrices $\bFb$, using Lemma \ref{Lem.5.1}. Then, in the next section, for a,  \emph{now smooth,} matrix $\bFb$,   we can use the pseudodifferential calculus  to establish the \emph{asymptotic estimate} for singular numbers of $\g\Yb^*,$ containing the integral norm of $\bFb$ but depending now only on the ellipticity constant and some algebraic norm of the approximating matrix $\ah$ and its  inverse, or, what is the same, on these bounds for the initial matrix $\aF$. Finally, we collect our estimates and establish the asymptotic singular numbers bound for $\bFb=\aF-\ah.$

 The important point here that  we have already got rid of the non-smoothness of  the coefficients of the operator, the latter stayed behind in the operator $\Xb,$ therefore our reasoning uses essentially the smoothness of the coefficients $\ah(y).$

 As a preparation, we note that we can suppose that $\bFb=0$ outside an (arbitrarily small) neighborhood of  the boundary $\Si$ of the domain $\Om.$ This observation has  been  used many times in papers on our topic, including  \cite{Sandgren}, \cite{Agr06}, \cite{Suslina99}. A simple explanation is that if $\dist(\supp \bFb, \partial \Om)>0$ then the integral  kernel of $\tilde{\Tb}{}^{-1}$ is a smooth function on the set $\supp \bFb\times \Si$ (which is separated from the diagonal) and the  eigenvalues of the corresponding operator decay faster than any power,
 \begin{equation}\label{locality}
   \nb^{\sup}(\theta, \g(\bFb\nabla\tilde{\Tb}{}^{-1} )^*)=0, \, \mbox{for}\, \mbox{any}\,\theta>0, \, \supp{\bFb}\cap \Si=\varnothing.
 \end{equation}
 Therefore, we may suppose that  $\bFb$ is supported in a conveniently  small neighborhood of $\Si$ where it possesses the properties discussed in Section \ref{Sect.geom}.
 \begin{proposition}\label{Thm.Est}
  Let $\ah(y), y\in \Om,$ be a smooth elliptic matrix function and $\bFb(y)\in L_\infty(\Om).$ Suppose that $\bFb$ is $L_\infty(\Si)$- continuous at $\Si\times \{1\}$ and  $\bFb$  is zero outside some neighborhood of $\Si\times \{1\}$.
  Then
  \begin{equation}\label{Yb88}
    \nb^{\sup}(2d, \g\Yb[\bFb]^*)\le C(\ah)\|\bFb(.,1)\|_{2d+2}^{\frac{d}{2}}\|\bFb(.,1)\|_{\infty}^{\frac{ d}{2}}.
  \end{equation}
\end{proposition}

\begin{proof}
  We consider the operator $\Zb[\bFb]=\g\Yb^*$, where $\Yb=\Yb[\bFb]=\bFb\nabla\tilde{\Tb}^{-1}$.
Since  we are studying the singular numbers of  $\Zb[\bFb]:L_2(\Om)\to L_2(\Si)$, we can consider the adjoint operator $\Zb[\bFb]^*:L_2(\Si)\to L_2(\Om)$ instead. This operator admits a convenient representation, namely,

\begin{equation}\label{W0}
   \Zb[\bFb]^*=(\g\Yb^*)^* =\left[\g(\bFb\nabla\tilde{\Tb}^{-1})^*\right]^*=
   (\bFb\nabla) (\g\tilde{\Tb}^{-1})^*.
 \end{equation}
 In fact, starting with the identity

 \begin{equation*}
   \bFb\nabla \tilde{\Tb}^{-1}=(\bFb\nabla \tilde{\Tb}^{-\frac12})\tilde{\Tb}^{-\frac12},
 \end{equation*}
 we obtain
 \begin{equation*}
( \bFb\nabla \tilde{\Tb}^{-1})^{*}= \tilde{\Tb}^{-\frac12}(\bFb\nabla \tilde{\Tb}^{-\frac12})^*;
 \end{equation*}
 further, we have
 \begin{equation*}
  \g( \bFb\nabla \tilde{\Tb}^{-1})^{*}=\left(\g\tilde{\Tb}^{-\frac12}\right)\left(\bFb\nabla \tilde{\Tb}^{-\frac12}\right)^*.
 \end{equation*}
And now we take  adjoints,
 \begin{gather*}
   \Zb(\bFb)^*\equiv\left[\g( \bFb\nabla \tilde{\Tb}^{-1})^{*}\right]^*=\\\nonumber
   \left(\bFb\nabla \tilde{\Tb}^{-\frac12}\right)^{**}\left(\g\tilde{\Tb}^{-\frac12}\right)^*=\bFb\nabla \tilde{\Tb}^{-\frac12}(\g\tilde{\Tb}^{-\frac12})^*=\bFb\nabla (\g\tilde{\Tb}^{-1})^*,
 \end{gather*}
 since $\bFb$ is a bounded function.
 We consider now  the function $w(y)=(\g\tilde{\Tb}{}^{-1})^*g\in L_2(\Om)$ for $g\in L_2(\Si).$ By Lemma \ref{Lem2}, this function satisfies the second order elliptic  equation $\Lch w=0$ in $\Om.$ Its restriction to $\Si$ is $\g w=\g(\g\tilde{\Tb}^{-1})^*g$, and, since by Proposition \ref{Prop.Pseudodiff} the operator $\g(\g\tilde{\Tb}^{-1})^*$ is an order $-1$ pseudodifferential operator,  the inequality holds
 \begin{equation}\label{W3}
   \|\g w\|_{H^{1}(\Si)}\le C\|g\|_{L_2(\Si)}.
 \end{equation}

 By the elliptic regularity property, for the  solution $w$ of  the second order elliptic equation $\Lch w=0$ in $\Om$ with the Dirichlet boundary condition $\g w$,  the Sobolev spaces estimate holds:
 \begin{equation}\label{W4}
   \|w\|_{H^{s+\frac12}(\Om)}\asymp \|\g w\|_{H^{s}(\Si)},
 \end{equation}
  now for \emph{all} $s\in\R^1.$ Note here, that the constants concealed in the $'\asymp'$ symbol in \eqref{W4} depend on the particular value of $s$ and, what is important,  on the operator $\Lch$, and may deteriorate when $\Lch$ is changing, while the derivatives of the coefficients of $\Lch$ grow in the process of approximation, even with the ellipticity constant preserved.

  Now we pass to estimating the singular numbers of the operator $(\g \Yb^*)^*$, whose \emph{squares}, by \eqref{W0}, are described by the  ratio
 \begin{equation}\label{W5}
   \frac{\int_{\Om}|\bFb(y)\nabla
   (\g\tilde{\Tb}^{-1})^*g|^2dy}{\|g\|^2_{L_2(\Si)}}=\frac{\int_{\Om}|\bFb(y)\nabla w(y)|^2dy}{\|g\|^2_{L_2(\Si)}},
 \end{equation}
 where, recall, $w=(\g \tilde{\Tb}^{-1})^* g$, $g\in L_2(\Si).$

 We use  the relations \eqref{W4} (with $s=1$) and \eqref{W3}. As a result,  the eigenvalues of the ratio \eqref{W5} are majorized by the eigenvalues of the ratio

 \begin{equation}\label{W6}
   \frac{\int_{\Om}|\bFb(y)\nabla w(y)|^2dy}{\|w\|^2_{H^{3/2}(\Om)}}, \, \Lch w=0\,  \mbox{{in}}\, \Om.
 \end{equation}

 On the next step, while evaluating the eigenvalues of the ratio \eqref{W6}, we use the following weighted estimate for solutions of elliptic equations, see \cite{Suslina99}, Lemma 3.3.
 \begin{lem}\label{Lem.Sus}
 Let a function $w\in H^s(\Om), \, s>0,$ be a solution of an elliptic equation $\Lch w=0$ with smooth coefficients in a bounded domain $\Om$ and let $r(y)$ be the distance from the point $y\in\Om$  to the boundary of $\Om.$ Then, for any number $\ka\ge 0$ such that $s+\ka$ is an integer,
 \begin{equation}\label{W7}
   \int_{\Om} r(y)^{2\ka} |\nabla_{s+\ka} w(y)|^2 dy\le C \|w\|^2_{H^s(\Om)},
 \end{equation}
 for a constant $C=C(s,\ka,\Om)$ not depending on $w.$
  \end{lem}
The constant in \eqref{W7} depends, of course, on the operator $\tilde{\Lc}$ as well, but this dependence is not mentioned in \cite{Suslina99} -- although this fact is of no importance at the moment, as long as  the operator is fixed.

We also have the obvious inequality,
\begin{equation*}
  \int_{\Om} r(y)^{2\ka} | \nabla w(y)|^2 dy\le C \|w\|^2_{H^s(\Om)},
\end{equation*}
for $s\ge 1.$

We consider the case $d>1$ first, thus excluding temporarily the case of a two-dimensional domain $\Om.$
We apply Lemma \ref{Lem.Sus} for the particular values $s=\frac32,$ $\ka=\frac12,$ $s+\ka=2,$ and thus reduce our task to  estimating the eigenvalues of the ratio

\begin{equation}\label{W9}
  \frac{\int_{\Om}|\bFb(y)\nabla w(y)|^2dy}{\int_{\Om} r(y) (|\nabla_2 w|^2+|\nabla w|^2 dy }, \Lch w=0;
\end{equation}
where, in our case, $r(y)=1-y_{d+1}$
If  we drop the condition $\Lch w=0$ in \eqref{W9}, the eigenvalue counting function may only increase and the result will give us an upper eigenvalue estimate for  the ratio \eqref{W9}.

The denominator  in \eqref{W9} is the quadratic form of  a degenerate elliptic operator. Such kind of spectral problems was considered in the series of papers of M.Z. Solomyak and I.L. Vulis in 1970-s, see, especially,  \cite{VuSolIzv},  where order sharp eigenvalue estimates and formulas for asymptotics have been proved. We need only estimates; we might have cited the paper \cite{Suslina99}, where the eigenvalue estimates for problems of the type \eqref{W9} were established in an even more general setting, namely, for a domain with  \emph{piecewise smooth} boundary, see there Lemma 4.1. However, for our, more simple, case, since  the boundary is smooth,  we refer to a more easily accessible and more elementary paper  \cite{VuSolIzv}. In order to do this, we make some more transformations of our spectral problem. Namely, we majorize the matrix $\bFb(y)$ by its matrix norm $b(y)=|\bFb(y)|=(\sum|b_{j,k}(y)|^2)^{\frac12}.$ After this, we replace the gradient of $w$, $\nabla w=(\partial_1w,\dots,\partial_{d+1}w)$  by an arbitrary vector function with $d+1$ components $\wb_\io$. This widens the set of functions where the variational ratio is considered, therefore, we arrive at the ratio
\begin{equation}\label{W10}
  \frac{\sum_\io \int_{\Om} b(y)^2 |\wb_\io(y)|^{2} dy}{\sum_\io \int_{\Om} r(y) (|\nabla \wb_\io(y)|^2+|\wb_\io|^2) dy}.
\end{equation}

In this way, the spectral problem \eqref{W9} splits into the direct sum of $d+1$ identical scalar spectral problems,
\begin{equation}\label{W11}
  \frac{ \int_{\Om} b(y)^2 |\wb(y)|^{2} dy}{ \int_{\Om} r(y) (|\nabla \wb(y)|^2+|\wb|^2) dy},
\end{equation}
exactly of the form, considered in \cite{VuSolIzv}.
The result  which we are going to use here is the combination of  Lemma 5.1 and Lemma 5.4  in \cite{VuSolIzv}.

We cite and further discuss these results, tailored for our special case. Recall that the boundary $\Si$ corresponds to $y_{d+1}=1, $ $r(y)=1-y_{d+1},$ and the corresponding change is made in the formulation.

Set $m=d+1.$ Let $\bb[\wb]$  be the quadratic form in the cylinder $\Cs\subset\Om\subset \R^{m}$:
\begin{equation}\label{VS1}
  \bb[\wb]=\int_{\Cs}r(y)^{\be}b(y)^2\wb(y)|^2 dy,
\end{equation}
and
\begin{equation}\label{VS2}
  \ab^{(\a)}[\wb]=\int_{\Om}r(y)^{\a}(|\nabla \wb|^2+|\wb|^2) dy.
\end{equation}
\begin{lem}\label{5.1}
\textbf{(Lemma 5.1, a), in \cite{VuSolIzv}.)}
  Let the coefficient $b(y',y_{m}),$ $1-\de\le y_m\le 1,$ not depend on $y_{m}$, namely, $b(y',y_m)=\check{b}(y'), \, y'\in\Si, \check{b}\in L_\infty(\Si).$
  Suppose that $\kb:=2-\a+\be> \frac{2}{m}$ (this case is called the \emph{strong degeneration} case in \cite{VuSolIzv}). Set $ \theta=\frac{m-1}{\kb}(=\frac{d}{\kb})$.
Then for the spectrum of the ratio
\begin{equation}\label{VS3}
  \frac{\bb[\wb]}{\ab^{(\a)}[\wb]}, \wb\in H^1(\Om),
\end{equation}
the estimate holds
\begin{equation}\label{VS4}
  \nb^{\sup}(\theta,\ref{VS3})\le C \|\check{b} \|_{L_{2m}(\Si)}^{\frac{m-1}{2}}\|\check{b}\|_{L_\infty(\Si)}^{\theta- \frac{m-1}{2}}.
\end{equation}
\end{lem}
The crucial importance of this estimate is that it involves the integral norm of $\check{b}.$

The proof of this Lemma in \cite{VuSolIzv} is based upon the separation of variables in the cylinder, thorough bookkeeping of the eigenvalues of the separated problem and finally using the result on the asymptotics of eigenvalues of an elliptic boundary problem with singular weight.

Note that in  \cite{VuSolIzv}, the weight function in the form $\bb$ is denoted by $b(y)$, while it is $b(y)^2$ in \eqref{VS1}; the corresponding change is made in the formulation of Lemma \ref{5.1}.

The second lemma relaxes the condition $b(y',y_{m})=\check{b}(y')$ used in Lemma \ref{5.1}. We formulate it in our terms:
\begin{lem}\label{5.4}\textbf{(Lemma 5.4 in \cite{VuSolIzv}.)}
  Let all conditions of Lemma \ref{5.1}, except $b=\check{b}(y')$,  be fulfilled, the latter being replaced by
  \begin{equation}\label{VS5}
   b(y',y_m)\, \mbox{is}\, L_\infty(\Si) \, \mbox{-continuous}\, \mbox{at}\, \Si\times\{1\}.
  \end{equation}
  Then the asymptotic estimate \eqref{VS4} holds, with $\check{b}=b(.,1)$.
\end{lem}
In fact, a  stronger statement is formulated in Lemma 5.4 in \cite{VuSolIzv}, namely, that the asymptotic bounds for $\nb_{\pm}$ for the functions $b$ and $\check{b}$ coincide, but  we need only the upper estimate \eqref{VS4}.

In \cite{VuSolIzv}, only a short sketch of the proof of Lemma 5.4 is given, with reference to 'standard tools of the variational method' (which might, in fact,  have been considered standard at that glory period of the variational method, but are, probably, not that standard nowadays). In more detail, and in much more generality, the reasoning, explaining the passage from $y_m$-independent coefficient $b$  to the one satisfying \eqref{VS5}, is given in \cite{Suslina99}, however, formally, the condition $b\in C(\Om)$ was imposed there. In fact,  only the condition \eqref{VS5} was actually used in \cite{Suslina99}.  A short but, hopefully, sufficient explanation is the following.

 By the Ky Fan inequality, it suffices to prove that if \eqref{VS5} holds and $b(.,1) =0$  then $\nb^{\sup}(\theta,\eqref{VS3})=0$.  Fix some $\pmb{\e}>0$ and find $l>0$ such that $|b(y)|<\pmb{\e}$ for $y_m\ge 1-l.$ The quadratic form ${\bb[\wb]}$  splits into the sum of two forms,
\begin{gather}\label{VS6}
  \bb[\wb]=\int_{y_m\ge 1-l} b(y)^2|\wb(y)|^2dy+\int_{y_m<1-l}b(y)^2|\wb(y)|^2dy\le \\\nonumber
  C\left(\pmb{\e}^2\int_{\Om}  |\wb(y)|^2dy +  l^{-\de}\int_{\Om}|y_m-1|^{\de}|\wb(y)|^2dy\right)\equiv\bb_{\pmb{\e}}[\wb]+\bb_{l}[\wb], \, \de>0.
\end{gather}

For the first term in the second line in \eqref{VS6}, Lemma \ref{5.1} with $\check{b}=\pmb{\e}$ applies by monotonicity, which gives the estimate
$\nb_\pm^{\sup}(\theta, \bb_{\pmb{\e}}/\ab)\le C \pmb{\e}^{2\theta}.$
 For the eigenvalues of the operator described by the second term,  $\bb_{l}[\wb],$ we apply   Lemma \ref{5.1} with some $\be=\de>0,$  this means, with a different order in the weight in the numerator \eqref{VS1}. As a result, Lemma \ref{5.1} gives for this term a faster eigenvalues decay: \eqref{VS4} takes the form $\nb_\pm^{\sup}(\theta', \bb_{\de,L}/\ab )<\infty,$ $\theta'=\frac{m-1}{1+\de}<\theta,$ and therefore, $\nb(\theta, \bb_{\de,l}/\ab )=0.$ After this, the required equality  $\nb_{\pm}^{\sup}(\theta,\eqref{VS3})=0$ follows due to the arbitrariness of $\pmb{\e}.$

We apply Lemma \ref{5.4} for  $m=d+1$, $\a=1$, $\be=0,$ $\kb=2-\a+\be=1$, $\theta =d$, to the ratio \eqref{W11} and this gives us the desired estimate for the singular numbers of $\g\Yb^*$

The above reasoning breaks down in the two-dimensional case, $m=d+1=2,$ since here  $\kb=1$ and we have the equality $\kb=\frac{2}{m}$ instead of the required inequality $\kb>\frac{2}{m}.$ In this case, called in \cite{VuSolIzv} 'the intermediate degeneration', the above scheme produces a non-sharp order in the singular numbers estimate. This kind of complication was handled in \cite{Suslina99}, Sect.4, in the following way which we adapt to our situation. The idea is in using \eqref{W7} for a different, larger, value of $\ka$, so that the resulting spectral problem becomes the one with strong degeneration. The order of eigenvalue estimates obtained in this way does not depend on the chosen value of $\ka$. (Note that the dimension of the enveloping space, denoted by $m$ in \cite{VuSolIzv}, is denoted by $m+1$ in \cite{Suslina99}, therefore we change notations correspondingly when citing the latter paper.)

For $m=2, d=1,$ we choose the number $\ka$ in \eqref{W7} to be  not $\frac12$ but $\frac32,$ so, $s+\ka=3$, and the weighted inequality  \eqref{W7} takes the form
\begin{equation}\label{W12}
  \int_{\Om}r(y)^{3}|\nabla_3 w|^2 d y\le C \|w\|^2_{H^{\frac32}(\Om)}, \, \Lch w=0.
\end{equation}
Following the reasoning above, we arrive at estimating the singular numbers of the ratio
\begin{equation}\label{W13}
  \frac{\int_\Om  |\bFb(y)|^2|\nabla w|^2 dy}{\int_{\Om}r(y)^3 (|\nabla_3 w(y)|^2+|w(y)|^2)dy}.
\end{equation}
instead of \eqref{W11}.
With this set of parameters, the problem \eqref{W13} in dimension $m=2$ is of a strong degeneration type, The result in \cite{Vul76A} gives in this case the estimate, see also \cite{Suslina99}, Lemma 4.1 or \cite{Suslina21}.
\begin{equation}\label{W14}
  \nb^{\sup}(1, \ref{W13})\le C(\aF) \|\bFb(.,1)\|_{L_{4}(\Si)}^{1/3}\|\bFb(.,1)\|_{L_\infty(\Si)}^{2/3}.
\end{equation}

Formally, the proof of Lemma 4.1 in \cite{Suslina99} requires $\bFb$ to be continuous, however the passage to the discontinuous $\bFb$ which is $L_\infty(\Si)$- continuous at $\Si$ as  function of $y_{d+1}$ variable is made identically with the above case of the dimension $d>1.$  A more simple treatment of this case, for a smooth boundary, can be found in  \cite{Vul76A}.

Finally, we recall that the eigenvalues of the ratio \eqref{W5} are squares of the singular numbers of the operator $\Zb[\bFb],$ therefore, the spectral estimate of order $d$ for \eqref{W5} produces the singular numbers estimate of order $2d$ for $\g\Yb^*$
\end{proof}

\section{Sharp estimates for $\Zb=\g\Yb^*$}
We recall  that the constant $C(\ah)$ in \eqref{Yb88} depends in a non-controllable way on the approximating  matrix  $\ah$ (although, in an analogous situation, in \cite{Suslina99} it was found that the constants in this kind of eigenvalue estimates depend only on the bounds for some finite collection of derivatives of $\ah.$)
 Now we, for the case of a \emph{smooth} matrix $\bFb(y),$ establish an estimate for $\nb^{\sup}(2d,\g\Yb^*)$ in terms of the norm of $\bFb$ in an integral metric and the bounds for  the  principal symbol of the operator $\Lch.$

 In this section we establish asymptotic bounds for singular numbers of the operator $\Zb$ with smooth matrix $\bFb\in C^{\infty}(\overline{\Om})$, which depend only on the ellipticity constant of the matrix $\ah$ (in other words, on the ellipticity constant of the operator $\Lc$), its norm and certain integral norm of the matrix $\bFb(.,1)$.

\subsection{The structure of $\Zb^*\Zb$}
The squares of  singular numbers of $\Zb$ are eigenvalues of the self-adjoint operator $\Zb^{*}\Zb$ acting in $L_2(\Si)$. We discuss its structure in more detail now.
We have
\begin{gather}\label{Yb1}
\Zb^*\Zb=\left[(\bFb(y)\nabla) (\g \tilde{\Tb}{}^{-1})^*\right]^*\left[\bFb(y)\nabla (\g \tilde{\Tb}{}^{-1})^*\right]\\\nonumber
= \g \tilde{\Tb}{}^{-1}[(\bFb(y)\nabla)^*(\bFb(y)\nabla)](\g \tilde{\Tb}{}^{-1})^*.
\end{gather}

Since all coefficient functions entering in $\Zb$ are smooth,  operators composing \eqref{Yb1} belong to the Boutet-de-Monvel
 algebra of pseudodifferential operators (see, e.g.,  \cite{Grubb.book}). We consider separate terms more closely.

The operator $\g \tilde{\Tb}{}^{-1}$ and its adjoint were described in detail in Sect. 4.4 and 4.5.
Since $\Rb$ is an integral operator with Hermitian kernel $R(y,z)$ the operator in \eqref{Yb8} is an integral operator  acting from $\Om$ to $\Si$ with kernel
\begin{equation*}
  R^{(1)}(y,z)=\g_y R(y,z)-\ND_y\g_{1,y}R(y,z), \, y\in\Si,\, z\in\Om,
\end{equation*}
where $\ND_y$, $\g_y,$ $\g_{1,y}$ denote the N-to-D operator $\ND,$ the trace $\g$ and the conormal derivative $\g_1$ acting upon the $y$ variable at $y\in \Si,$ $\g_{1,y}=\g_y \partial_{\n_{\ah}(y)}.$ Recall that $\ND$ is an order $-1$ pseudodifferential operator at the boundary with symbol $\be(y',\y'),$ see \eqref{beta}.

The  adjoint operator $(\g \tilde{\Tb}{}^{-1})^*$ is, therefore, an integral  operator acting from $\Si$ to $\Om$ with  the adjoint kernel,

\begin{equation*}\label{Yb94}
R^{(1)*}(y,z)=R^{(1)}(z,y)=\g_zR(y,z)-\ND_z\g_{1,z}R(y,z), \, y\in \Om, z\in\Si.
\end{equation*}
Now we can describe how the operator $\Zb =\bFb(y)\nabla (\g \tilde{\Tb}{}^{-1})^*$ acts. By \eqref{Yb94}, $\Zb =\bFb(y)\nabla (\g \tilde{\Tb}{}^{-1})^*$ is an integral operator acting from $\Si$ to $\Om$ with integral kernel
   \begin{equation}\label{Yb10}
     \Zc(y,z)=\bFb(y)\nabla_y\g_zR(y,z)- \bFb(y)\nabla_y\ND_z\g_{1,z}R(y,z), \, y\in \Om, z\in \Si.
   \end{equation}
Similarly, $\Zb^*$ is an integral operator acting from $\Om$ to $\Si$ with the integral  kernel
\begin{equation}\label{Yb11}
  \Zc^*(y,z)=\Zc(z,y).
\end{equation}

\subsection{Composition}\label{Sect.Comp} Now we collect the description of the entries in the operator $\Wb=\Zb^*\Zb.$ By  \eqref{Yb10},\eqref{Yb11}, $\Zb^*\Zb$ is the composition
\begin{equation}\label{Yb12}
  \Zb^*\Zb=((\bFb\nabla)((\g-\ND\g_1)\Rb)^*)^*(\bFb\nabla)((\g-\ND\g_1)\Rb)^*
\end{equation}
is the integral operator with kernel
\begin{equation}\label{Yb13}
  \Wc(y,z)=\int_{\Om}\Zc(y,\varsigma)\Zc(\varsigma,z)d\varsigma,
\end{equation}
with kernel $\Zc$ given by \eqref{Yb10}, acting on $\Si.$

As elements in the Boutet-de-Montvel algebra, $(\g-\ND\g_1)\Rb$ is a trace operator acting from $\Om$ to $\Si,$  its adjoint
$((\g-\ND\g_1)\Rb)^*$ is a Poisson operator acting from $\Si$ to $\Om,$ and the whole composition, $\Wb,$ has integral kernel of the form
\begin{gather}\label{Yb14}
  \Wc(y,z)=\int_{\Om}\Zc(y,\sigma)\Zc(\sigma,z)d\sigma= \\\nonumber
  \int_{\Om} (\g_y-\ND_y (\g_1)_y)(\bFb(\sigma)\nabla_{\sigma})^*R(y,\sigma)\bFb(\sigma)\nabla_\sigma (\g_z-\ND_z(\g_1)_z)R(\sigma,z) d\sigma.
\end{gather}
 By the Boutet de Monvel calculus, $\Wb$  is a pseudodifferential operator of order $-1$  on $\Si.$

We need to find the dependence of the principal symbol of this operator on the matrices $\bFb(y)$ and $\aF(y),$ $y\in\Sigma.$ In this evaluation, we may ignore terms of lower order, in particular, those appearing when we commute the factors in the product in \eqref{Yb14}.

We start by considering the term $(\bFb(\sigma)\nabla_\sigma)^*(\bFb(\sigma)\nabla_\sigma).$ This is a pseudodifferential operator of order $2$ in $\Om$ with principal symbol satisfying
\begin{equation}\label{b.nabla}
  \Xc(\sigma,\varsigma)=\langle\bFb(\sigma)\varsigma,\bFb(\sigma)\varsigma\rangle=|\bFb(\sigma)\varsigma|^2\le |\bFb(\sigma)|^2|\varsigma|^2
\end{equation}

Further on, the truncated self-adjoint pseudodifferential operator $\Mb_{\bFb}=$ \- $ (\bFb\nabla \Rb)^*(\bFb\nabla \Rb) $ of order $-2$ has principal symbol

\begin{equation}\label{b.nabla.R}
   \Mc_\bFb(\sigma,\varsigma)=\rF_{-2}(\sigma,\varsigma)^2 |\bFb(\sigma)\varsigma|^2,
\end{equation}
where, recall,  $\rF_{-2}(\sigma,\varsigma)=\langle\ah(\sigma)\varsigma,\varsigma\rangle^{-1}$ is the principal symbol of the fundamental solution $\Rb$, so $\Mb_{\bFb}$  is an operator of order $-2$.
Next, we need to make the restriction of $\Mb_{\bFb} $ to the boundary, $y\in\Si.$

The  trace operators $\g-\ND\g_1$ acting from both sides upon  $\Mb\equiv\Mb_{\bFb}$ produce an order $-1$ pseudodifferential operator on the boundary with principal symbol, again, containing the matrix  $\bFb^*\bFb$ and a homogeneous symbol depending algebraically  on the matrix $\ah$ and its inverse $\ah^{-1}.$ Its principal symbol can be calculated following the Boutet de Monvel calculus rules, which involve only algebraical operations with principal symbols. We however do not need to calculate this symbol or write down its explicit expression. For our needs it is sufficient to describe its properties, especially, on what data of our operator it depends. We demonstrate it on one of the terms,  $\Mb_{\bFb}^{(\g)}:=\g(\g \Mb_{\bFb})^*$, which is the restriction of $\Mb_{\bFb}$ to the boundary. Its symbol  can be evaluated as

\begin{gather}\label{RbbRgamma}
  \Mc^\g(y,\y')=\frac{1}{2\pi}\int_{-\infty}^\infty\Mc_{\bFb}(y,\y',\y_{d+1})d\y_{d+1}\\\nonumber
  =\frac{1}{2\pi}\int_{-\infty}^\infty  |\bFb(y)\y|^2(\ah(y)\y,\y)^{-4} d\y_{d+1} \le \\\nonumber
 |\bFb(y)|^2 |\ah(y)|^{-2}\frac{1}{2\pi}\int_{-\infty}^\infty (|\y'|^2+|\y_{d+1}|^2)^{-1}d \y_{d+1}\le C|\ah(y)^{-2}| |\bFb(y)|^2|\y'|^{-1},
  \, y\in \Si,
\end{gather}
 the expression \eqref{RbbRgamma} is calculated
in local co-ordinates at the point $y\in \Si $, where $y_{d+1}$ axis is directed along the normal to $\Si.$

In a similar way, other terms in the symbol of $\Zb^*\Zb$ can be estimated, using the representation \eqref{Yb12}, following the composition rules in the Boutet de  Monvel algebra, see, e.g.,  \cite{Grubb.book}, Sect.9.5, 10.4.  Here we use the expression for the symbol of the conormal derivative, $\imath\langle\ah(y)\y,\n\rangle$ and the symbol $\beta(y,\y')$ of the $\ND$ operator, see \eqref{beta}. All these terms are majorated by $|\bFb(y)|^2,$ with  bounded dependence on $|\ah(y)^{-1}|$ and $|\ah(y)|,$ $y\in\Si.$ Singular Green operators arising in the process  of composition of truncated pseudodifferential operators give no contribution to the principal symbol of the composition

 As a result we have
\begin{proposition}\label{prop.princ.symb}
  The principal symbol $\WF(y,\y')$, $(y,\y')\in T^*\Si$ of the order $-1$  pseudodifferential operator $\Wb=\Zb^*\Zb$  admits the estimate
  \begin{equation}\label{princ.symbol}
   |\WF(y,\y')|\le C |\bFb(y)|^2 A(|\ah(y)^{-1}|,|\ah(y)|)|\y'|^{-1},
  \end{equation}
 with a function $A(s_1,s_2)$ bounded on bounded intervals separated from zero.
\end{proposition}

Now we apply the asymptotic estimate for singular numbers of negative order pseudodifferential operators. This formula was obtained by M.Birman and M.Solomyak in \cite{BSpseudo}, Theorem 1, for operators in a domain in the Euclidean space and then carried over to manifolds in \cite{BSSib}. The proof in  \cite{BSpseudo}, II, a very technical one, having been published in the, now quite obscure, Russian journal, was almost unaccessible to Western researchers, although the result  was rather widely cited. Fortunately, quite recently, R.Ponge in \cite{Ponge}, Sect.6, proposed a rather soft proof of a special case of the main results in \cite{BSpseudo}, concerning operators with smooth symbols, which fits our needs.
We arrive at the following asymptotic estimate for the singular numbers  of the operator $\Zb.$

\begin{proposition}\label{Prop.Yb}
  Let $\tilde{\Lc}$ be an elliptic operator with leading coefficients matrix $\ah(y)\in C^\infty(\overline{\Om}).$ Suppose that for all $y\in \overline{\Om},$
  \begin{equation}\label{unif.Estim.L}
    |\ah(y)^{-1}|, |\ah(y)|\le \Cb_{\aF}.
  \end{equation}
  Let $\bFb(y)$ be a smooth symmetric matrix in $C^\infty(\overline{\Om}).$
  Then for the singular numbers of the operator $\Zb=\bFb\nabla(\g\tilde{T}^{-1})^*$ the estimate holds
  \begin{gather}\label{Eigenv.estim.Yb}
    \nb^{\sup}(2d,\Zb)=\nb^{\sup}(d, \Zb^*\Zb)\le A(\Cb_{\aF}) \int_{\Si}|\bFb(y')|^{2d}d\m_\Si(y')\le   \\\nonumber
    C  A(\Cb_{\aF})\|\bFb\|_{L_{2d+2}(\Si)}^{d}\|\bFb\|_{L_\infty(\Si)}^{d},
  \end{gather}
  with constant $A(\Cb_{\aF})$ depending only on the  $L_{\infty}$ norms of $\ah^{-1}$ and $\ah$.

\end{proposition}
\begin{proof} As we just found, the operator $\Zb^*\Zb$ is an order $-1$ pseudodifferential operator on the boundary $\Si,$ with leading symbol $\WF(y',\y'),$ $(y',\y')\in \mathrm{T}^*\Si$ majorated by $\Cb_{\aF}|\bFb(y)|^2|\y'|^{-1}.$ By the formula (23) in \cite{BSpseudo}, for the eigenvalue counting function $n(\la,\Zb^*\Zb)$ the asymptotics holds (in our notations)
\begin{equation}\label{BSAsympt}
  \nb(2d,\Zb)=\nb(d,\Zb^*\Zb)=d^{-1}(2\pi)^{-d}\int_\Si\int_{\mathrm{S}^*\Si}\WF(y,\y')^{d}d\m_{\Si}\om(\y'),
  \end{equation}
  and $\om$ is the standard volume form on the unit sphere $S^{d-1}.$
Thus, our estimate for the principal symbol $\WF$ of the operator $\Zb^*\Zb$, substituted in \eqref{BSAsympt}, gives \eqref{Eigenv.estim.Yb}. The last inequality in \eqref{Eigenv.estim.Yb} follows from the fact that the $L_\infty$ norm of $\bFb$ majorates its $L_p$ norm, for any $p<\infty.$
\end{proof}

\subsection{Final estimates for $\g\Yb^*$}
We combine the results of the last two subsections in order to obtain the final singular numbers estimate for $\Zb=\g \Yb^*.$

\begin{proposition}\label{Prop.Limit} Let $\bFb(y)$ be a matrix function satisfying the conditions of Proposition \ref{Thm.Est} and $\ah(y)$ be a smooth matrix satisfying $|\ah^{-1}(y)|, |\ah^(y)|\le \Cb_{\ah}.$
Then for the singular numbers of  the operator $\g(\Yb[\bFb])^*=\g(\bFb\nabla(\tilde{\Tb}^{-1}))^*$ the singular numbers estimate holds
\begin{equation}\label{lim.good}
  \nb(2d,\g\Yb[\bFb]^*)\le A( \Cb_{\ah})\|\bFb(.,1)\|_{L_{2d+2}(\Si)}^{{d}}.
\end{equation}
\end{proposition}
Note that, from the first glance, the statement of Propostion \ref{Prop.Limit} coincides with the one of Proposition \ref{Thm.Est}. There is, however a critical improvement. While the constant in \eqref{Yb88} may depend, in an uncontrollable way, on the matrix $\ah(y),$  the constant in \eqref{lim.good} depends only on the ellipticity constant $|\aF^{-1}|$ and the norm  $|\ah(y)|,$ but not on the matrix $\ah$ itself or its derivatives.

\begin{proof} For a given $\e,$  we construct, following the procedure in Lemma \ref{lem.appr}, the approximating \emph{smooth} matrix $\bFb_\e\in C^{\infty}(\overline{\Om})$ such that $\|\bFb_\e\|_{L_\infty(\Si)} \le C_0  \|\bFb_\e\|_{L_\infty(\Om)}$ and
\begin{equation}\label{lim.good.1}
  \|\bFb(.,1)-\bFb_\e(.,1)\|_{L_{2d+2}(\Si)}<\e.
\end{equation}
  By the Ky Fan inequality,
  $n(\la_1+\la_2, \Kb_1+\Kb_2)\le n(\la_1,\Kb_1)+n(\la_2,\Kb_2),$
 it  follows in the usual way (see, e.g. \cite{BS}, Sect.5), that
  \begin{equation}\label{lim.good.2}
    \nb^{\sup}(2d,\g\Yb[\bFb]^*)^{(1+2d)^{-1}}\le \nb^{\sup}(2d,\g\Yb[\bFb_\e]^*)^{(1+2d)^{-1}}+ \nb^{\sup}(2d, \g\Yb[\bFb-\bFb_\e]^*)^{(1+2d)^{-1}}.
  \end{equation}
  The  last term in  \eqref{lim.good.2} tends to zero as $\e\to 0$ by Proposition \ref{Thm.Est}, applied to $\bFb-\bFb_\e$  in place of $\bFb$.  Therefore, passing to the limit in \eqref{lim.good.2}, we obtain
  \begin{equation}\label{lim.good.3}
\nb^{\sup}(2d,\g\Yb[\bFb]^*)\le     \limsup_{\e\to 0}\nb^{\sup}(2d,\g\Yb[\bFb_\e]^*).
  \end{equation}
  Now we  remember the estimate \eqref{Eigenv.estim.Yb}, which is valid, since $\bFb_\e$ is smooth now. By the construction of the approximating matrix $\bFb_\e$ in Lemma \ref{lem.appr}, the norm $\|\bFb_\e(.,1)\|_{L_{2d+2}(\Si)}$ is controlled by the same norm  of $\bFb$. Therefore we can pass to limit as $\e\to 0$ in the inequality
  \begin{equation}\label{lim.good.4}
    \nb^{\sup}(2d,\g\Yb[\bFb]^*)\le
    C\lim\sup_{\e\to 0}\|\bFb_\e(.,1)\|^{d}_{L_{2d+2}(\Si)}\le C\|\bFb(.,1)\|^{\frac{d}{2}}_{L_{2d+2}(\Si)},
  \end{equation}
  which gives the required estimate, since the norm $\|\bFb(.,1)\|_{L_\infty}\le C\|\aF\|_{L_\infty}$ can be absorbed in the coefficient $C.$
\end{proof}

\section{Conclusion of the proof}

The uniformity property of the estimate in \eqref{lim.good} enables us to establish the crucial approximation result for the operators $\Hb_{\ah}.$

\begin{thm}\label{Approximation.Thm} Let $\aF$ be a matrix satisfying the conditions of  Lemma \ref{lem.appr}, and for a given $\e>0$,  $\ah=\ah_\e$ be the approximating smooth matrix constructed in this Lemma. Then
\begin{equation}\label{Conv.thm.1}
  \nb^{\sup}(2d,\g(\Yb[\aF-\ah])^*)\le C(\aF) \e^{d},
\end{equation}
and
\begin{equation}\label{Conv.thm.2}
 \lim_{\e\to 0} \nb^{\sup}(d, \Hb_{\aF}-\Hb_{\ah})\to 0.
\end{equation}
\end{thm}
\begin{proof}
We recall the representation $\Hb_{\aF}-\Hb_{\ah}=(\g(\Yb[\aF-\ah])^*)^*(\g\Xb^*).$
  Therefore, the second estimate will follow in the usual way, compare, e.g.,  \cite{BS}, from the first one, by the Ky Fan inequality for the product of operators,
 \begin{equation}\label{KyFanProd}
  n(\la_1\la_2, \Kb_1\Kb_2)\le n(\la_1,\Kb_1)+n(\la_2,\Kb_2),
 \end{equation}
 where  $\Kb_1=\g\Xb^*, \, \Kb_2=(\g(\Yb[\aF-\ah])^*)^*,$
   since for $\Kb_1$
 we already know the estimate $\nb^{\sup}(2d,\g\Xb^*)<\infty$.
 Namely, for a given $\la>0,$ we set $\la_1=\e^{-d/4}\la,$ $\la_2=\e^{d/4}\la$ in \eqref{KyFanProd}, $\la^2=\la_1\la_2,$ then  both terms on the right in \eqref{KyFanProd} get a  small factor as $\e\to 0.$

 To prove \eqref{Conv.thm.1}, we apply Proposition \ref{Prop.Limit} with $\bFb=\aF-\ah.$
\end{proof}
Finally,
using the basic asymptotic perturbation lemma, Lemma \ref{Lem.5.1}, we establish our main result.

\begin{proof} Of Theorem \ref{RTthm}. Let $\aF$ be the coefficient matrix in our P-S problem in a smooth domain. Consider its smooth approximation  $\ah$  constructed according to Lemma \ref{lem.appr}, with $p=2d+2$. For the operator $\Hb_{\ah}$ the asymptotic formula \eqref{formula} is known. For the difference, $\Hb_\aF-\Hb_{\ah}$ we have the estimate \eqref{Conv.thm.2}. In these conditions, Lemma \ref{Lem.5.1} grants that for the limit operator $\Hb_{\aF}$ the eigenvalue asymptotics is valid and the coefficient in the asymptotics is given by the limit in the formula \eqref{formula}, using Lemma 2.5.
\end{proof}
\appendix
\section{Potential theory approach}We consider here the Poincar\'{e}-Steklov problem for the Laplacian, first, in dimension $d+1\ge 3.$ Consider the single and double layer potential operators on the Lipschitz boundary   $\Si,$
\begin{equation}\label{SL}
  \Ss: L_2(\Si)\to L_2(\Si), \, \Ss: f(x)\mapsto\int_\Si R(x-y)f(y) d\m_\Si(y),
\end{equation}
and
\begin{equation}\label{DL}
  \Ds: L_2(\Si)\to L_2(\Si), \, \Ds:  f(x)\mapsto\int_\Si \partial_{\n(y)}R(x-y)d\m_{\Si}(y),
\end{equation}
where $R$ is the fundamental solution for the Laplacian in $\R^{d+1}.$
The N-to-D operator $\ND$ is expressed via these operators as
\begin{equation}\label{NDpot}
  \ND=\Ss(\frac12+\Ds)^{-1},
\end{equation}
see, e.g., \cite{Agr.UMN}, \cite{AgrTMMO}, \cite{Grubb.book}.
To be more exact, the $\ND$ operator is considered here on functions, orthogonal in $L_2(\Si)$ to constants. For a smooth surface $\Si,$ both potential operators are order $-1$ pseudodifferential operators on $\Si.$ Being considered on functions orthogonal to constants, the operator $(\frac12+\Ds)$ is invertible and

\begin{equation}\label{NDinv}
  (\frac12+\Ds)^{-1}=2-2\Ds(\frac12+\Ds)^{-1}.
\end{equation}
This means that
\begin{equation}\label{ND2}
  \ND=2\Ss - 2 \Ss \Ds(\frac12+\Ds)^{-1},
\end{equation}
Therefore the difference between the N-to-D operator $\ND$ and twice the single layer potential $2\Ss$ is an order $-2$ pseudodifferential operator, and, by standard perturbational arguments, the eigenvalue asymptotics of $\ND$ is the same in the leading term as the asymptotics for the single layer potential $2\Ss$; finding the latter is a simple exercise for a smooth surface.

If $\Si$ is not infinitely smooth, for example, belongs to $C^{1+\a}$, $\a>0$ (such surfaces are often called Lyapunov ones), the above scheme still works. First of all, discarding the pseudodifferential approach, one should return to considering the operators $\Ss$ and $\Ds$ as integral operators with kernel having singularity at the diagonal. Here, on the one hand, the eigenvalue asymptotics of the operator $\Ss$ is known since the results of M.Sh. Birman and M.Z. Solomyak on  general weakly polar integral operators, see \cite{BSint}, where such eigenvalue asymptotics was found under rather weak regularity conditions. At the same time, the operator $\Ds$ is now a singular integral operator, but for a surface of class $C^{1+\a}$, it is still compact. Thus, the operator $\ND$ becomes, by \eqref{ND2}, a relatively compactly perturbed
single layer operator, and, again, by standard methods (say, using the Keldysh perturbation theorem) we arrive at the Weyl asymptotics for the $\ND$ operator.

So, it looks like this is not that long way remains to Lipschitz surfaces. However some serious obstacles appear. First of all, for a Lipschitz surface, the double layer integral operator $\Ds$ is not obviously bounded. In fact, it required a hard work to prove that it is bounded in $L_2(\Si)$ and, moreover, that the operator $\frac12+\Ds$ is invertible. Further on, it turned out that $\Ds$ is not necessarily compact; even in the two-dimensional case, in the presence of a corner point of a piece-wise smooth boundary, it has a fragment of essential spectrum near zero. As for the single layer potential $\Ss,$ the results of \cite{BSint} are, by themselves, not sufficient to establish eigenvalue asymptotics. However, in \cite{AgrAmos}, M.S.Agranovich and B.A.Amosov succeeded in proving order sharp two-sided eigenvalue estimates for potential type operators
on Lipschitz surfaces, which, in their turn, produce two-sided eigenvalue estimates for $\ND.$ Moreover, a localization technic was developed in \cite{AgrAmos}. It  turned out that if the non-smooth singularities of the surface $\Si$ are localized in the sense that the surface is smooth outside a closed set of zero measure, then the contribution of these singularities of $\Si$ both to the behavior of eigenvalues of $\Ss$ and to the non-compactness of $\Ds$ is negligible, and the Weyl asymptotic formula for the P-S problem was proved in this setting.

In  that period, after M.S.Agranovich started popularizing the problem on Steklov eigenvalue asymptotics for Lipschitz surfaces without additional restrictions, Prof. Grigory Tashchiyan and the author of this paper  managed in 2006 to dispose of one of complications. Namely, in \cite{RT1}, we used the approximation of Lipschitz surfaces by smooth ones, in the sense of a certain integral norm, to prove that the corresponding potential type integral operators, properly relocated, so that they act on one and the same surface, can be compared,  their difference becomes small in the spectral sense (see Lemma \ref{5.1}), and the passage to limit in formulas for spectral asymptotics becomes possible. Although this result did not solve the Lipschitz Steklov eigenvalue asymptotics problem, it solved a closely related one. Namely, in the \emph{interface} spectral problem, where the spectral parameter is placed, for solutions of the (Laplace, in the simplest case) equation at the \emph{jump} of the normal derivative across the Lipschitz surface $S\subset \Om:$
\begin{gather}\label{jump}
  \Om=\Om_1\cup\Om_2\cup S, \, S=\overline{\Om_1}\cap\overline{\Om_2}, \, \Delta u=0,\, \mbox{in}\, \Om_1\cup\Om_2,\\\nonumber
  \lambda\left[ \frac{\partial u}{\partial_{\n(x)}}\right]_{x}=u(x), \, x\in S,
\end{gather}
the solution involves only the single layer potential operator $\Ss$ and the Weyl asymptotics for the Lipschitz surface is justified.   Rather optimistically, we hoped that the Steklov problem will succumb as well.

One of the approaches we tried to explore was the following. Although we know that the double layer operator $\Ds$ is not expected to be compact, probably, the \emph{composition } $\Ds\Ss$ can be shown to be \emph{weaker} than $\Ss$ itself. The notion 'weaker' should imply that the singular numbers of $\Ds\Ss$ decay faster than the ones of $\Ss$,
\begin{equation}\label{minor}
  s_n(\Ds\Ss)=o(s_n(\Ss))=o(n^{-\frac1d}).
\end{equation}
  In the eigenvalue studies, there are a lot of occasions when the multiplication by a non-compact operator  nevertheless improves the rate of decay of the eigenvalues of the product.

\textbf{Hypothesis}. Consider the  composition $\Ds\Ss$ of single and double layer potentials. Recall that $\Ss$ is the integral operator with kernel $C|x-y|^{1-d}$ over a Lipschitz surface of dimension $d>1.$ Let $\Tc(x,y)$ be the integral kernel of the composition $\Ss\Ds$,
\begin{equation}\label{composit}
  \Tc(x,y)=\int_{\Si}|x-z|^{1-d}\partial_{\n(z)}|z-y|^{1-d} dz.
\end{equation}
Suppose, hopefully, that
\begin{equation}\label{Comp.Hyp}
  \Tc(x,y)|x-y|^{d-1}=o(1) \,\mbox{as }\, y\to x, \, \mbox{uniformly }\, \mbox{in}\, x\in \Si.
\end{equation}

\begin{proposition}\label{Prop.Hypot}
  Suppose that $d\ge 3$ and  \eqref{Comp.Hyp} is satisfied. Then $\la_j(\Ss\Ds)=o(\la_j(\Ss)).$
\end{proposition}
\begin{proof} For a fixed $\e$, we find a neighborhood $\UF_\e\subset\Si\times\Si$ of the diagonal $x=y$ such that $|\Tc(x,y)|<\e|x-y|^{1-d}$ for $(x,y)\in \UF_\e$. Take a
  cut-off function $\chi_\e(x,y)$ which equals zero outside $\UF_\e$ and equals $1$ inside, this means, in a neighborhood of the diagonal $x=y.$ The operator $\Ts=\Ss\Ds$ splits into the sum of two operators, $\Ts_\e,$ with kernel $\Tc(x,y)\chi_\e(x,y)$ and $\Ts_\e'$ with kernel $\Tc(x,y)(1-\chi_\e(x,y)).$  The second operator has a bounded kernel, therefore, in particular, belongs to the Hilbert-Schmidt class, and its singular values decay at least as $\la_j(\Ts_\e')=O(j^{-\frac12}), $ i.e. faster than the eigenvalues of $\Ss.$ As for the first operator, $\Ts_\e,$ we can apply the result by G.Kostometov, \cite{Kostometov}, Theorem 1 (see also, \cite{AgrAmos}, Theorem 4.1.) According to this theorem, for an integral operator $\Ts$ in a bounded domain in $\R^d$ with  kernel having  form $\Tc(x,y)=\phi(x,y)|x-y|^{-\kb},$ $d/2 < \kb<d,$ with a bounded function $\phi$, the eigenvalue estimate holds
  \begin{equation}\label{Kost.Est}
    \nb^{\sup}(\theta, \Ts)\le C \|\phi\|_{L_{\infty}}^{\theta}, \theta =\frac{d}{d-\kb.}
  \end{equation}
  In our case, in local co-ordinates, for $\Tc_\e(x,y)=\Tc(x,y)\chi_\e(x,y),$ $\kb=d-1$ and $\phi(x,y)=\Tc(x,y)\chi(x,y)|x-y|^{1-d}$, $\|\phi\|_{L_{\infty}}\le \e$. So, \eqref{Kost.Est} gives $\nb^{\sup}(\theta, \Ts_\e)\le C \e^{\theta},$ and by the arbitrariness of $\e,$ we arrive at the required estimate.
\end{proof}This estimate leads to the justification of the Weyl asymptotics, provided \eqref{Comp.Hyp}.

Using somewhat finer estimates, one covers the case of $d=2$ (excluded in the estimates in \cite{Kostometov}, \cite{AgrAmos}, since here the condition $d/2<\kb$ is not met). The hypothesis itself, meanwhile, stays unresolved, although for some particular nonregular surfaces we were able to justify it.

\section{Relaxing the smoothness conditions.}
Another possible approach consists in finding the conditions on the surface $\Si$, just a little bit more restrictive than Lipschitz ones, but still granting the compactness of the double layer operator $\Ds$, restricted to the complement of some small set, containing major singularities.

Here, the results of \cite{HMT1} proved to be useful.  Let $\Om\subset\R^{d+1}.$ The domain $\Om$ is called $\mbox{VMO}_1$-domain if \eqref{Lip2} holds  and the almost everywhere existing  normal vector field $\pmb{\n}(x)$ belongs to the space $\mbox{VMO}(\Si)$, the closure of the space of continuous vector fields on $\Si$ in the  $\mbox{VMO}$ metric. Of course, the class of $\mbox{VMO}_1$ domains contains all $C^1$ domains, but does not contain \emph{all} Lipschitz domains.
By Theorem 4.35 in \cite{HMT1}, singular integral operators of the type of double layer potential $\Ds,$ are compact in $L_2(\Si).$ By a localization of this property, as it was done, e.g., in  (see, e.g., \cite{AgrAmos}, \cite{AgrTMMO}, \cite{Agr06}) we arrive at the following result.
\begin{proposition} Let the $\Om$ be a Lipschitz domain and let  the boundary $\Si$ belong  to $\mbox{VMO}_1$ outside a closed set $E$ of the surface measure zero. Then for the $\ND$  operator the Weyl asymptotic formula holds.
\end{proposition}
The proof, leading, again, to \eqref{minor}, consists in splitting the Steklov problem $\partial_\n u=\la^{-1} u$ on $\Si$ into two weighted ones,
\begin{equation}\label{splitvmo}
  \partial_\n u=\la \ro_\e u,\,\mbox{and}\,  \partial_\n u=\la(1- \ro_\e) u,
\end{equation}
where $\ro_e$ is the characteristic function of a small neighborhood of the set $E$, having surface measure less than $\e$. This leads to the corresponding splitting of the operator $\Hb$ in \eqref{Permuted} into $\Hb=\Hb_{\ro_\e}+\Hb_{1-\ro_\e} .$ For the first of these operators, the smallness for $\nb^{\sup}(d,\Hb_{\ro_\e})$ can be proved on the base of Theorem \ref{Th.estimate}. For the remaining operator, $\Hb_{1-\ro_\e},$ it follows from the $VMO_1$ compactness property for the double layer potential that the $\nb^{\sup}(d,\Hb_{\ro_\e})=0.$

One may hope that the approximation of a Lipschitz surface by $VMO_1$ ones may lead to a more straightforward proof of the Weyl asymptotics for the P-S eigenvalues.

\end{document}